\documentclass[10pt]{article}
\usepackage{amsmath,amssymb,amsthm,amscd}
\numberwithin{equation}{section}

 \DeclareMathOperator{\Aut}{Aut}

\def\Bl{\biggl(}
\def\Br{\biggr)}
\def\p{\partial}
\def\b{\bar}

\def\o{\omega}
\def\l{\lambda}

\def\Tilde{\widetilde}

\def\tr{\rm tr}

\def\RR{{\mathbb R}}
\def\CC{{\mathbb C}}

\newtheorem{prop}{Proposition}[section]
\newtheorem{theo}[prop]{Theorem}
\newtheorem{lem}[prop]{Lemma}
\newtheorem{cor}[prop]{Corollary}
\newtheorem{rem}[prop]{Remark}

\newtheorem{defi}[prop]{Definition}
\newtheorem{conj}[prop]{Conjecture}
\newtheorem{q}[prop]{Question}
\def\begeq{\begin{equation}}
\def\endeq{\end{equation}}

\def\and{\quad{\rm and}\quad}

\let\lra=\longrightarrow

\def\mapright\#1{\,\smash{\mathop{\lra}\limits^{\#1}}\,}

\begin {document}
\bibliographystyle{plain}
\title{Ricci flow on K\"ahler-Einstein manifolds}
\author{X.X. Chen and G. Tian}
\date{Aug. 26, 2000 \\Revised on May 20, 2001}
\maketitle
\tableofcontents
\section{Introduction}

This is the continuation of our earlier paper \cite{chentian001}.
For any K\"ahler-Einstein surfaces with positive scalar curvature,
if the initial metric has positive bisectional curvature, then we
proved \cite{chentian001} that the K\"ahler-Ricci flow converges
exponentially to a unique K\"ahler-Einstein metric in the end.
This answers partially to a long standing problem in Ricci flow:
 on a compact K\"ahler-Einstein manifold, does the
K\"ahler Ricci flow converge to a K\"ahler-Einstein metric if the
initial metric has positive bisectional curvature? In this paper,
we will give a complete affirmative answer to this problem.

\begin{theo}  Let $M$ be a K\"ahler-Einstein manifold with
positive scalar curvature. If the initial metric has nonnegative
bisectional curvature and positive at least at one point, then
K\"ahler Ricci flow will converge exponentially fast to a
K\"ahler-Einstein metric with constant bisectional curvature.
\end{theo}
\begin{rem}
This problem was completely solved by R. Hamilton in the case of
Riemann surfaces (cf. \cite{Hamilton88}).  We also refer the
reader to B. Chow's paper \cite{Chow91} for more developments on
this problem.
\end{rem}
As a direct consequence, we have the following:
\begin{cor}
The space of K\"ahler metrics with non-negative bisectional
curvature (and positive at least at one point) is path-connected.
The space of metrics with non-negative curvature operator (and
positive at least at one point)  is also path-connected.
\end{cor}


\begin{theo}  Let $M$ be any K\"ahler-Einstein orbifold (cf. Definition 9.2)  with positive
scalar curvature. If the initial metric has non-negative
bisectional curvature and positive at least at one point, then
the K\"ahler Ricci flow converges exponentially to a
K\"ahler-Einstein metric with constant bisectional curvature.
Moreover, $M$ is a global quotient of $\CC P^n.\;$
\end{theo}

Clearly, Corollary 1.3 holds in the case of K\"ahler orbifolds.

\begin{rem}
What we really need is that the Ricci curvature is positive along
the K\"ahler-Ricci flow. Since the positivity on Ricci curvature
may not be preserved under the Ricci flow, we will use the fact
that the positivity of the bisectional curvature is preserved.
\end{rem}

\begin{rem} In view of the solution of the Frankel
conjecture solved by S. Mori \cite{Mori79} and Siu-Yau
\cite{Siuy80}, it suffices to study this problem on a K\"ahler
manifold which is biholomorphic to $\CC P^n$. However, we don't
need to use the result of Frankel conjecture.  Moreover,  we do
not use explicitly the knowledge of the positive bisectional
curvature. We use this condition only when we quote a result of
Mok and Bando (cf. \cite{chentian001}), and a classification
theorem by M. Berger.
\end{rem}

\begin{rem} We need the assumption on the existence of K\"ahler-Einstein
metric because we will use a nonlinear inequality from
\cite{tian98}. Such an inequality is just the
Moser-Trudinger-Onofri inequality if the underlying  manifold is
the Riemann sphere.
\end{rem}

\begin{rem}  If we assume the existence of a lower bound for the
functional $E_1 - E_0$\footnote{cf. Section 2.3 for definition of
$E_0, \;E_1.$ }, then we shall be able to derive a convergence
result similarly. Therefore, it is interesting to study the lower
bound of $E_1-E_0$ among metrics whose bisectional curvature is
positive.
\end{rem}

\begin{rem} We learn from H. D. Cao \cite{caoprivate} that the
holomorphic orthogonal bisectional curvature \footnote{It is the
bisectional curvature between two any two orthogonal complex plan.
}is preserved under the K\"ahler Ricci flow (this will follow from
Mok's proof by a simple modification).  It is easy to see that
positive Ricci curvature is preserved under the flow. Then our
proof will extend to this case. Note that the bisectional
curvature is not necessary positive during the flow.
\end{rem}

 Now let us review briefly the history of Ricci flow.  The
Ricci flow was first introduced by R. Hamilton in
\cite{Hamilton82}, and it has been a subject of intense study ever
since. The Ricci flow provides an indispensable tool of deforming
Riemannian metrics towards to canonical metrics, such as Einstein
ones. It is hoped that by deforming a metric to a canonical
metric, one can  understand the geometric and topological
structures of underlying manifolds. For instance, it was proved
\cite{Hamilton82} that any closed 3-manifold of positive Ricci
curvature is diffeomorphic to a spherical space form. We refer the
reader to \cite{Hamilton93} for more information.

If the underlying manifold is a K\"ahler manifold, the Ricci flow
preserves the K\"ahler class.
 Following a similar idea of Yau \cite{Yau78}, Cao \cite{Cao85} proved
that the solution converges to a K\"ahler-Einstein metric if the
first Chern class of the underlying K\"ahler manifold is zero or
negative. Consequently, he re-proved the famous Calabi-Yau
theorem\cite{Yau78}. On the other hand, if the first Chern class
of the underlying K\"ahler manifold is positive, the solution of
the K\"ahler Ricci flow may not converge to any K\"ahler-Einstein
metric. This is because there are compact K\"ahler manifolds with
positive first Chern class which do not admit any
K\"ahler-Einstein metrics (cf. \cite{futaki83}, \cite{Tian97}). A
natural and challenging problem is whether or not the K\"ahler
Ricci flow on a compact K\"ahler-Einstein manifold converges to a
K\"ahler-Einstein metric. Our theorem settles this problem in the
case of K\"ahler metrics of positive bisectional curvature or
positive curvature operator. It was proved by S. Bando
\cite{Bando84} for 3-dimensional K\"ahler manifolds and by N. Mok
\cite{Mok88} for higher dimensional K\"ahler manifolds that the
positivity of bisectional curvature is preserved under the
K\"ahler Ricci flow.

The typical method in studying the Ricci flow depends on
pointwise bounds of the curvature tensor by using its evolution
equation as well as the blow-up analysis. In order to prevent
formation of singularities, one blows up the solution of the Ricci
flow to obtain profiles of singular solutions. Those profiles
involve Ricci solitons and possibly more complicated singular
models. Then one tries to exclude formation of singularities by
checking that these solitons or models do not exist under
appropriate global geometric conditions. It is a common sense
that it is very difficult to detect how the global geometry
effects those singular models even for a very simple manifold like
$\CC P^2$. The first step is to classify those singular models
and hope to find their geometric information. Of course, it is
already a very big task. There have been many exciting works on
these (cf. \cite{Hamilton93}).

Our new contribution is to find a set of new functionals which
are the Lagrangians of certain new curvature equations involving
various symmetric functions of the Ricci curvature. We show that
these functionals decrease essentially along the K\"ahler Ricci
flow and have uniform lower bound. By computing their
derivatives, we can obtain certain integral bounds on curvature
of metrics along the flow.

For the reader's convenience, we will recall what we study in
\cite{chentian001} regarding these new functionals. In
\cite{chentian001},  we proved that the derivative of each $E_k$
along an orbit of automorphisms gives rise to a holomorphic
invariant $\Im_k$, including the well-known Futaki invariant as a
special one. When $M$ admits a K\"ahler-Einstein metric, all
these invariants $\Im_k$ vanish, and the functionals $E_k$ are
invariant under the action of automorphisms.

Next we  proved in \cite{chentian001}  that these $E_k$ are
bounded from below.  We then computed the derivatives of $E_k$
along the K\"ahler-Ricci flow. Recall that the K\"ahler Ricci
flow is given by
\begin{equation}
   {{\partial \varphi} \over {\partial t }} =  \log {(\omega + \partial \overline\partial
\varphi)^n \over {\omega}^n } + \varphi - h_{\omega},
\end{equation}
where $h_\omega$ depends only on $\omega$. The derivatives of
these functionals are all bounded uniformly from above along the
K\"ahler Ricci flow. Furthermore, we found that $E_0$ and $E_1$
decrease along the K\"ahler Ricci flow. These play a very
important role in this and the preceding paper. We can derive from
these properties of $E_k$ integral bounds on curvature, e.g. for
almost all K\"ahler metrics $\omega_{\varphi(t)}$ along the flow,
we have

\begin{equation}
 \int_M\; (R(\omega_{\varphi(t)})-r)^2 \;{\omega_{\varphi(t)}}^n \rightarrow 0,
\label{eq:keyestimate2}
\end{equation}
where $R({\omega_{\varphi(t)}})$ denotes the scalar curvature and
$r$ is the average scalar curvature.

In complex dimension 2,  using
 the above integral bounds on the curvature with Cao's Harnack inequality and the
generalization of Klingenberg's estimate, we can bound the
curvature uniformly along the K\"ahler Ricci flow in the case of
K\"ahler-Einstein surfaces. However, it is not enough in high
dimension, since the formula (\ref{eq:keyestimate2}) is not
scaling invariant.  We must find a new way of utilizing this
inequality in higher dimensional manifolds.  Following the work of
C. Sprouse \cite{spro001},  J. Cheeger  and T. Colding \cite{CC96}
of deriving a uniform upper-bound on the diameter,  we then use a
result of Li-Yau \cite{liyau80} and a theorem of C. Croke
\cite{croke80} to derive a uniform upper bound on both the Sobolev
constant and the Poincare constant on the evolved K\"ahler metric.
Once these two important constants are bounded uniformly,  we can
use the  Moser iteration to obtain $C^0$ estimate along the
modified K\"ahler Ricci flow.
A priori, this curve of evolved K\"ahler-Einstein metrics is not
even differentiable on the level of potentials in terms of time
parameter. This gives us a lot of troubles in deriving the desired
$C^0$ estimates.  What we need is  to re-adjust this curve of
automorphisms so that it is at least $C^1$ uniform on the level of
K\"ahler potentials. Once $C^0$ estimate is established,  it is
then possible to obtain the $C^2$ estimate (following a similar
calculation of Yau \cite{Yau78}) and  Calabi's $C^3$ estimates.
Eventually, we can prove that the modified K\"ahler Ricci flow
converges exponentially to the unique K\"ahler-Einstein metric.

Unlike \cite{chentian001}, we don't use any pointwise estimate on
curvature; in particular, we don't need to use the Harnack
inequality. It appears to us that the fact that the set of
functionals we found being essentially decreasing along the
K\"ahler Ricci flow and  having a uniform lower bound at the same
time, has already exclude the possibilities of formation of
singularities. In higher dimensional manifolds, this idea of
having integral estimates on curvature terms,
may prove to be an effective and attractive alternative (vs. the usual  pointwise estimates). \\

In this paper, we also extend our results to K\"ahler-Einstein
orbifolds with positive bisectional curvature. Note that the limit
metric of the K\"ahler-Ricci flow on orbifolds must be Einstein
metric with positive bisectional curvature. M. Berger's theorem
\cite{Berger65} then implies that it must be of constant
bisectional curvature. We then use the exponential map to
explicitly prove that such an orbifold must be
a global quotient of $\CC P^n.\;$\\

 The organization of our paper is roughly as follows: In Section 2, we
review briefly some basics in K\"ahler geometry and some results
we obtained in \cite{chentian001}.  In Section 3, we prove that
for any K\"ahler metric in the canonical class with non-negative
Ricci curvature,  if the scalar curvature is sufficiently closed
to the average in the $L^2$ sense, then it has uniform diameter
bound. Next using the old results  of Li-Yau and the result of C.
Croke, we bound both the Sobolev constant and the Poincar\'e
constant. In Section 4, we prove $C^0$ estimates for  all time
over the
modified K\"ahler Ricci flow. 
 In Section 5, we prove that we can choose a uniform
gauge. In Section 6, we obtain both $C^2$ and $C^3$ estimates. In
section 7, we prove the exponential convergence to the unique
K\"ahler-Einstein metric with constant bisectional curvature.
 In Section 8,  we prove that any orbifold supports a K\"ahler metric with positive constant bisectional curvature
 is globally a quotient of $\CC P^n.\;$ In Section 9, we prove
 Theorem 1.5 and make some concluding remarks and propose some open questions.
\section{Setup and known results}
\subsection{Setup of notations}

Let $M$ be an $n$-dimensional compact K\"ahler manifold. A K\"ahler metric can be given by its
K\"ahler form $\omega$ on $M$. In local coordinates $z_1, \cdots, z_n$, this $\omega$ is of
the form
\[
\omega = \sqrt{-1} \displaystyle \sum_{i,j=1}^n\;g_{i
\overline{j}} d\,z^i\wedge d\,z^{\overline{j}}  > 0,
\]
where $\{g_{i\overline {j}}\}$ is a positive definite Hermitian matrix function.
The K\"ahler condition requires that $\omega$ is a closed positive
(1,1)-form. In other words, the following holds
\[
 {{\partial g_{i \overline{k}}} \over
{\partial z^{j}}} =  {{\partial g_{j \overline{k}}} \over
{\partial z^{i}}}\qquad {\rm and}\qquad {{\partial g_{k
\overline{i}}} \over {\partial z^{\overline{j}}}} = {{\partial
g_{k \overline{j}}} \over {\partial
z^{\overline{i}}}}\qquad\forall\;i,j,k=1,2,\cdots, n.
\]
The K\"ahler metric corresponding to $\omega$ is given by
\[
 \sqrt{-1} \;\displaystyle \sum_1^n \; {g}_{\alpha \overline{\beta}} \;
d\,z^{\alpha}\;\otimes d\, z^{ \overline{\beta}}.
\]
For simplicity, in the following, we will often denote by $\omega$ the corresponding K\"ahler metric.
The K\"ahler class of $\omega$ is its cohomology class $[\omega]$ in $H^2(M,\RR).\;$
By the Hodge theorem, any other K\"ahler
metric in the same K\"ahler class is of the form
\[
\omega_{\varphi} = \omega + \sqrt{-1} \displaystyle \sum_{i,j=1}^n\;
{{\partial^2 \varphi}\over {\partial z^i \partial z^{\overline{j}}}}
> 0
\]
for some real valued function $\varphi$ on $M.\;$ The functional
space in which we are interested (often referred as the space of
K\"ahler potentials) is
\[
{\cal P}(M,\omega) = \{ \varphi \;\mid\; \omega_{\varphi} = \omega
+ \sqrt{-1}
 {\partial} \overline{\partial} \varphi > 0\;\;{\rm on}\; M\}.
\]
Given a K\"ahler metric $\omega$, its volume form  is
\[
  \omega^n = {1\over {n!}}\;\left(\sqrt{-1} \right)^n \det\left(g_{i \overline{j}}\right)
 d\,z^1 \wedge d\,z^{\overline{1}}\wedge \cdots \wedge d\,z^n \wedge d\,z^{\overline{n}}.
\]
Its Christoffel symbols are given by
\[
  \Gamma^k_{i\,j} = \displaystyle \sum_{l=1}^n\;g^{k\overline{l}} {{\partial g_{i \overline{l}}} \over
{\partial z^{j}}} ~~~{\rm and}~~~ \Gamma^{\overline{k}}_{\overline{i}\,\overline{j}} =
\displaystyle \sum_{l=1}^n\;g^{\overline{k}l} {{\partial g_{l \overline{i}}} \over
{\partial z^{\overline{j}}}}, \qquad\forall\;i,j,k=1,2,\cdots n.
\]
The curvature tensor is
\[
 R_{i \overline{j} k \overline{l}} = - {{\partial^2 g_{i \overline{j}}} \over
{\partial z^{k} \partial z^{\overline{l}}}} + \displaystyle \sum_{p,q=1}^n g^{p\overline{q}}
{{\partial g_{i \overline{q}}} \over
{\partial z^{k}}}  {{\partial g_{p \overline{j}}} \over
{\partial z^{\overline{l}}}}, \qquad\forall\;i,j,k,l=1,2,\cdots n.
\]
We say that $\omega$ is of nonnegative bisectional curvature if
\[
 R_{i \overline{j} k \overline{l}} v^i v^{\overline{j}} w^k w^{\overline{l}}\geq 0
\]
for all non-zero vectors $v$ and $w$ in the holomorphic tangent
bundle of $M$. The bisectional curvature and the curvature tensor
can be mutually determined. The Ricci curvature of $\omega$ is
locally given by
\[
  R_{i \overline{j}} = - {{\partial}^2 \log \det (g_{k \overline{l}}) \over
{\partial z_i \partial \bar z_j }} .
\]
So its Ricci curvature form is
\[
  {\rm Ric}(\omega) = \sqrt{-1} \displaystyle \sum_{i,j=1}^n \;R_{i \overline{j}}(\omega)
d\,z^i\wedge d\,z^{\overline{j}} = -\sqrt{-1} \partial \overline{\partial} \log \;\det (g_{k \overline{l}}).
\]
It is a real, closed (1,1)-form. Recall that $[\omega]$ is called
a canonical K\"ahler class if this Ricci form is cohomologous to
$\lambda \;\omega,\; $ for some constant $\lambda.$

 \subsection{The K\"ahler Ricci flow}

    Now we assume that the first Chern class $c_1(M)$ is positive.
The normalized Ricci flow (c.f. \cite{Hamilton82} and
\cite{Hamilton86}) on a K\"ahler manifold $M$ is of the form
\begin{equation}
  {{\partial g_{i \overline{j}}} \over {\partial t }} = g_{i \overline{j}}
  - R_{i \overline{j}}, \qquad\forall\; i,\; j= 1,2,\cdots ,n,
\label{eq:kahlerricciflow}
\end{equation}
if we choose the initial K\"ahler metric $\omega$ with $c_1(M)$ as
its K\"ahler class. The flow (2.1) preserves the K\"ahler class
$[\omega]$. It follows that on the level of K\"ahler potentials,
the Ricci flow becomes
\begin{equation}
   {{\partial \varphi} \over {\partial t }} =  \log {{\omega_{\varphi}}^n \over {\omega}^n } + \varphi - h_{\omega} ,
\label{eq:flowpotential}
\end{equation}
where $h_{\omega}$ is defined by
\[
  {\rm Ric}(\omega)- \omega = \sqrt{-1} \partial \overline{\partial} h_{\omega}, \; {\rm and}\;\displaystyle \int_M\;
  (e^{h_{\omega}} - 1)  {\omega}^n = 0.
\]
As usual, the flow (2.2) is referred as the K\"ahler Ricci flow
on $M$.

The following theorem was proved by S. Bando for 3-dimensional
compact K\"ahler manifolds. This  was later proved by N. Mok in
\cite{Mok88} for all dimenisonal K\"ahler manifolds. Their proofs
used Hamilton's maximum principle for tensors. The proof for
higher dimensions is quite intrigue.

\begin{theo} \cite{Bando84}  \cite{Mok88} Under the K\"ahler Ricci flow, if the initial metric
has nonnegative bisectional curvature, then the evolved metrics
also have non-negative bisectional curvature. Furthermore, if the
bisectional curvature of the initial metric is positive at least
at one point, then the evolved metric has positive bisectional
curvature at all points.
\end{theo}

Before Bando and Mok, R. Hamilton proved (by using his maximum
principle for tensors)

\begin{theo} Under the Ricci flow, if the initial metric
has nonnegative curvature operator, then the evolved metrics also
has non-negative curvature operator. Furthermore, if the
curvature operator of the initial metric is positive at least at
one point, then the evolved metric has positive curvature operator
at all points.
\end{theo}

\subsection{Results from the previous paper \cite{chentian001}}
In this subsection, we collect a few  results in our earlier paper
\cite{chentian001}. First, we introduce the new functionals $E_k =
E_k^0 - J_k (k =0,1,2\cdots, n) $ where $E_k^0$ and $J_k$ are
defined below.

\begin{defi} For any $k=0,1,\cdots, n$, we define a functional $E_k^0$
on ${\cal P}(M,\omega)$ by
\[
  E_{k,\omega}^0 (\varphi) = {1\over V}\; \displaystyle \int_M\;  \left( \log {{\omega_{\varphi}}^n \over \omega^n}
   - h_{\omega}\right) \left(\displaystyle \sum_{i=0}^k\; {{\rm Ric}(\omega_{\varphi})}^{i}\wedge\omega^{k-i} \right)
    \wedge {\omega_{\varphi}}^{n-k} + c_k,
\]
where
\[
c_k ={1\over V}\; \displaystyle \int_M\;
   h_{\omega} \left(\displaystyle \sum_{i=0}^k\; {{\rm Ric}(\omega)}^{i}\wedge\omega^{k-i} \right)
    \wedge {\omega}^{n-k}.
\]
\end{defi}

\begin{defi}
For each $k=0,1,2,\cdots, n-1$, we will define $J_{k, \omega}$ as
follows: Let $\varphi(t) $ ($t\in [0,1]$) be a path from $0$ to
$\varphi$ in ${\cal P}(M,\omega)$, we define
\[ J_{k,\o}(\varphi) = -{n-k\over V}
\int_0^1 \int_M {{\partial \varphi}\over{\partial t}}
\left({\omega_{\varphi}}^{k+1} - {\omega}^{k+1}\right)\wedge
{\omega_{\varphi}}^{n-k-1}\wedge dt.
\]
Put $J_n =0 $ for convenience in notations.
\end{defi}

\begin{rem}
In a non canonical K\"ahler class, we need to modify the
definition slightly since $h_{\omega}$ is not defined. For any
$k=0,1,\cdots, n,\;$ we define
\[
\begin{array}{lcl}
  E_{k,\omega}(\varphi) & = &
 {1\over V}\; \displaystyle \int_M\;   \log {{\omega_{\varphi}}^n \over \omega^n}
\; \left(\displaystyle \sum_{i=0}^k\; {{\rm
Ric}(\omega_{\varphi})}^{i}\wedge {\rm Ric}(\omega)^{k-i} \right)
    \wedge {\omega_{\varphi}}^{n-k}  \\
    &  & \qquad \qquad  -
    {{n-k}\over V} \displaystyle \int_M\; \varphi
    \left({\rm Ric}(\omega)^{k+1} - \omega^{k+1}\right) \wedge \omega^{n-k-1}
    - J_{k,\omega}(\varphi).
    \end{array}
\]
The second integral on the right hand side is to offset the change
from $\omega$ to $Ric(\omega)$ in the first term. The derivative
of this functional is exactly the same as in the canonical
K\"ahler class. In other words, the Euler-Lagrange equation is
not changed.
\end{rem}
If $\omega \in c_1(M),\;$ then we assume $E_k = E_{k,\omega}.\;$
Direct computations lead to
\begin{theo}For any $k=0,1,\cdots, n$, we have
\begin{eqnarray}
{d E_k \over dt} & = & {{k+1}\over V} \displaystyle \int_M
\Delta_{\varphi}\left( {{\partial \varphi}\over {\partial t}}
\right )\; {{\rm Ric}(\omega_{\varphi})}^{k} \wedge
{\omega_{\varphi}}^{n-k}
\nonumber \\
& & \qquad -
 {{n-k}\over V}\displaystyle \int_M {{\partial \varphi}\over {\partial t}} \left({{\rm Ric}(\omega_{\varphi})}^{k+1}
 - {\omega_{\varphi}}^{k+1}\right) \wedge  {\omega_{\varphi}}^{n-k-1}.
\label{eq:decay functional0}
\end{eqnarray}
Here $\{\varphi(t)\}$ is any path in ${\cal P}(M,\omega)$.
\end{theo}

\begin{prop} Along the K\"ahler Ricci flow, we have
\begin{equation}
{d E_k \over d t} \leq - {{k+1}\over V} \displaystyle \int_M
(R(\omega_\varphi)-r) {\rm Ric}(\omega_{\varphi})^{k} \wedge
{\omega_{\varphi}}^{n-k}. \label{eq: Ekdecreases}
\end{equation}
When $k=0 , 1$, we have
\begin{eqnarray}
{{d E_0 }\over{d\,t}} & =  & -{{n\sqrt{-1}}\over V} \displaystyle
\int_M \partial {{\partial \varphi}\over {\partial t}} \wedge
\overline{\partial} {{\partial \varphi}\over{\partial
t}}  {\omega_{\varphi}}^{n-1}\leq 0,\\
{{d E_1 }\over{d t}} & \leq & - {{2}\over V} \displaystyle \int_M
(R(\omega_\varphi)-r)^2 {\omega_{\varphi}}^{n} \leq  0. \nonumber
\end{eqnarray}
In particular, both $E_0$ and $E_1$ are decreasing along the
K\"ahler Ricci flow.
\end{prop}

 We then prove that the
derivatives of these functionals along a holomorphic automorphisms
give rise to holomorphic invariants. For any holomorphic vector
field $X,$ and for any K\"ahler metric $\omega,$ there exists a
potential function $\theta_X$ such that
\[
   L_{X} \omega = \sqrt{-1} \p \bar \p \theta_X.
\]
Here $L_X$ denotes the Lie derivative along a vector field $X$ and
$\theta_X$ is defined up to the addition of any constant. Now we
define $\Im_k(X,\omega)$ for each $k=0,1,\cdots, n$ by
\[
\begin{array}{ll}
\Im_k (X,\omega)  = (n-k) \displaystyle \int_M\theta_X\;
{\omega}^n\\
\qquad  + \displaystyle \int_M \left( (k+1) \Delta \theta_X
\;{{\rm Ric}(\omega)}^{k}\wedge {\omega}^{n-k} - (n-k)\; \theta_X
\;{{\rm Ric}(\omega)}^{k+1} \wedge {\omega}^{n-k-1}\right ).
\end{array}
\]
Here and in the following, $\Delta$ denotes the Laplacian of
$\omega$. Clearly, the integral is unchanged if we replace
$\theta_X$ by $\theta_X + c$ for any constant $c$.

The next theorem assures that the above integral gives rise to a
holomorphic invariant.

\begin{theo} The integral $\Im_k (X,\omega)$ is independent of choices of
K\"ahler metrics in the K\"ahler class $[\omega].\;$  That is,
$\Im_k (X,\omega)=\Im_k (X,\omega')$ so long as the K\"ahler forms
$\omega$ and $\omega'$ represent the same K\"ahler class. Hence,
the integral $\Im_k (X,\omega)$ is a holomorphic invariant, which
will be denoted by $\Im_k (X,[\omega])$.
\end{theo}

\begin{cor}
The above invariants $\Im_k (X, c_1(M))$ all vanish for any
holomorphic vector fields $X$ on a compact K\"ahler-Einstein
manifold. In particular, these invariants all vanish on $\CC P^n$.
\end{cor}

\begin{cor} For any K\"ahler Einstein manifold, $E_k (k=0,1,\cdots, n)$ is invariant under
actions of holomorphic automorphisms.
\end{cor}

One crucial step  in \cite{chentian001} is to modify the
K\"ahler-Einstein metric so that the evolved K\"ahler form is
centrally positioned with respect to this new K\"ahler Einstein
metric.  For the convenience of a reader, we include the
definition of ``centrally positioned" here.

\begin{defi}  Any K\"ahler form  $ \omega_{\varphi}$  is called centrally positioned
with respect to some K\"ahler-Einstein metric $\omega_{\rho} =
\omega + \sqrt{-1} \partial \overline{\partial} \rho$
 if it satisfies the following:
\begin{equation}
\displaystyle \int_M (\varphi -\rho) \; \theta\; {\omega_{\rho}}^n
= 0, \qquad \forall \;\theta \in \Lambda_1(\omega_{\rho}).
\label{eq:propposition}
\end{equation}
\end{defi}

\begin{prop} Let $\varphi(t)$ be the evolved K\"ahler potentials. For any $t> 0,$ there always exists
an automorphism $\sigma(t) \in {\rm Aut} (M)$ such that
$\omega_{\varphi(t)}$ is centrally positioned with respect to
$\omega_{\rho(t)}.\;$ Here
\[
  \sigma(t) ^* \omega_1 = \omega_{\rho(t)} = \omega + \sqrt{-1} \p
  \bar \p \rho(t),
\]
\end{prop}
where $\omega_1$ is an K\"ahler-Einstein metric.

\begin{rem} In \cite{chentian001}, we proved that the existence of at least one
K\"ahler-Einstein metric $\omega_{\rho(t)}$ such that
$\omega_{\varphi(t)}$ is centrally positioned with respect to
$\omega_{\rho(t)}.\;$ As a matter of fact, such a
K\"ahler-Einstein metric is unique.  However, a priori we don't
know if the curve $\rho(t)$ is differentiable or not.
\end{rem}
\begin{prop} On a K\"ahler-Einstein manifold, the $K$-energy $\nu_{\omega}$ is uniformly
bounded from above and below along the K\"ahler Ricci flow.
Moreover, there exists a uniform constant $C$ such that
\[
\begin{array}{ccl}
|J_{k,\omega_{\rho(t)}}(\varphi(t) - \rho(t) )| & \leq &
\{\nu_{\omega}(\varphi(t)) + C\}^{1\over \delta },\\
\log { {\omega_{\varphi}}^n \over {{\omega_{\rho(t)}}^n}} & \geq &
 - 4 C''\, e^{2\left(\nu_{\omega}(\varphi(t)) + C\right)^{1\over \delta} + C')},\\
E_k (\varphi(t)) & \geq &  - e^{c\left(1 + {\rm max} \{0,
\nu_\omega(\varphi(t))\} + \left(\nu_{\omega}(\varphi(t)) +
C\right)^{1\over \delta }\right) },
\end{array}
\]
where $c,\;C,\;C'$ and $C''$ are some uniform constants.  And
$\rho(t)$ is defined in the preceding proposition.
\end{prop}

\begin{cor} The energy functional $E_k (k=0,1,\cdots, n) $ has a
uniform lower bound from below along the K\"ahler Ricci flow.
 \end{cor}

\begin{cor} For each $k=0,1,\cdots, n,$ there exists a uniform constant $C$
such that the following holds (for any $T\leq \infty$) along the
K\"ahler Ricci flow:
\[
\displaystyle \int_{0}^{T}\; {{k+1}\over V}\; \displaystyle
\int_M\; \left(R(\omega_{\varphi(t)}) - r\right) \;  {{\rm
Ric}(\omega_{\varphi(t)})}^{k} \wedge {\omega_{\varphi(t)}}^{n-k}
\;d\,t \leq C.
\]
When $k=1, $ we have
\[
  \displaystyle \int_0^{\infty}\;{{1}\over V}\;
 \displaystyle \int_M\; (R(\omega_{\varphi(t)})-r) ^2\;   {\omega_{\varphi(t)}}^{n}  \;d\,t
 \leq C < \infty.
\]
\end{cor}

\section{Estimates of Sobolev and Poincare constants}
In this section, we will prove that for any K\"ahler metric in
the canonical K\"ahler class, if the scalar curvature is close
enough to a constant in $L^2$ sense and if the Ricci curvature is
non-negative, then there exists a uniform upper bound for both the
Poincar\'e constant and the Sobolev constant. We first follow an
approach taken by  C. Sprouse \cite{spro001} to obtain a uniform
upper bound on the diameter.

 In \cite{CC96}, J. Cheeger
and T. Colding proved an interesting and useful inequality which
converts integral estimates along geodesic to integral estimates
on the whole manifold.  In this section, we assume $m= {\rm dim}
(M).\;$
\begin{lem} \cite{CC96}
Let $A_1, \; A_2$ and $W$ be open subsets of $M$ such that $A_1,
A_2 \subset W,$ and all minimal geodesics $r_{x,y}$ from $x \in
A_1$ to $y \in A_2$ lie in $W.\;$ Let $f$ be any non-negative
function. Then
\[
\begin{array} {l} \displaystyle \int_{A_1 \times A_2} \displaystyle
\int_{r_{x,y}} \;  f(r(s)) \;d\,s\; d \; vol_{A_1 \times A_2}  \\
\qquad \leq C(m,k,\Re)({\rm diam}(A_2) vol(A_1) + {\rm diam}(A_1)
vol(A_2)) \displaystyle \int_W \; f\; d\;vol,
\end{array}
\]
where for $k \leq 0,$
\begin{equation}
C(m,k,\Re) = {{{\rm area} (\partial B_k (x,\Re))} \over {{\rm
area} (\partial B_k (x,{\Re \over 2}))}},
\end{equation}
\begin{equation}
\Re\geq \displaystyle \sup \{d(x,y) \mid (x,y) \in (A_1 \times
A_2)\},
\end{equation}
and $B_k(x,r)$ denotes the ball of radius $r$ in the simply
connected space of constant sectional curvature $k.\;$
\end{lem}
In this paper, we always assume ${\rm Ric} \geq 0$ on $M,$ and
thus $C(n,k,\Re) = C(n).\;$ \\
Using this theorem of Cheeger and Colding, C. Sprouse
\cite{spro001} proved an interesting lemma:

\begin{lem} \cite{spro001}
Let $(M,g)$ be a compact Riemannian manifold with ${\rm Ric} \geq
0.\;$ Then for any $\delta > 0$ there exists $\epsilon =
\epsilon(n, \delta) $ such that if
\begin{equation}
{1 \over V} \displaystyle \int_M \left((m-1) - Ric_- \right)_+ <
\epsilon(m,\delta), \label{eq:dimaeter0}
\end{equation}
then the $diam(M) < \pi + \delta.$ Here ${\rm Ric}_-$ denotes the
lowest eigenvalue of the Ricci tensor; For any function $f$ on
$M,\; f_+(x) = \max\{f(x),0\}.$
\end{lem}

\begin{rem} Note that the right hand side of equation (\ref{eq:dimaeter0}) is not
scaling correct. A scaling correct version of this lemma should
be: For any positive integer $a> 0,$ if
\[
{1\over V}  \displaystyle \int_M \; |Ric- a| \; d\,vol <
\epsilon(m,\delta) \cdot a,
\]
then the diameter has a uniform upper bound.
\end{rem}

\begin{rem} It is interesting to see what the optimal constant
$\epsilon(m,\delta)$ is.  Following this idea, the best constant
should be
\[
\epsilon(m,\delta)  = \displaystyle \sup_{N > 2}\; {{N-2}\over {8
C(m) N^m}}.
\]
However, it will be  interesting to figure out the best constant
here. 
\end{rem}

Adopting his arguments, we will prove the similar lemma,
\begin{lem} Let $(M,\omega)$ be a polarized K\"ahler manifold and
$[\omega]$ is the canonical K\"ahler class. Then there exists a
positive constant $\epsilon_0$ which only depends on the
dimension, such that if the Ricci curvature of $\omega$ is
non-negative and if
\[
 {1 \over V} \displaystyle \int_M \; (R -n)^2 \omega^n \leq
 \epsilon_0^2,
\]
then there exists a uniform upper bound on diameter of the
K\"ahler metric $\omega.\;$ Here $r$ is the average of the scalar
curvature.
\end{lem}
\begin{proof}
We first  prove that the Ricci form is close to its K\"ahler form
in the $L^1$ sense (after proper rescaling).  Note that
\[ {\rm Ric}(\omega) - \omega = \sqrt{-1} \partial \bar \partial f
\]
for some real valued function $f.\;$ Thus
\[
\displaystyle \int_M\; \left( {\rm Ric}(\omega) - \omega \right)^2
\wedge \omega^{n-2} = \displaystyle \int_M\; \left( \sqrt{-1}
\partial \bar \partial f\right)^2 \wedge \omega^{n-2} =0.
\]

On the other hand,  we have
\[ \displaystyle \int_M\;
\left( {\rm Ric}(\omega) - \omega \right)^2 \wedge \omega^{n-2} =
{1\over {n(n-1)}}\; \displaystyle \int_M\; \left((R - n)^2 - |{\rm
Ric}(\omega) - \omega |^2 \right) \omega^n.
\]
Here we already use the identity $tr_{\omega}  \left( {\rm
Ric}(\omega) - \omega \right) = R - n.\;$ Thus
\[
\displaystyle \int_M\;|{\rm Ric}(\omega) - \omega |^2  \; \omega^n
=\displaystyle \int_M\;(R - n )^2  \; \omega^n.
\]
This implies that
\[
\begin{array}{lcl} \left(\displaystyle \int_M |{\rm Ric} - 1
|\;\omega^n\right)^2 & \leq & \displaystyle \int_M |{\rm
Ric}(\omega) -
\omega |^2 \;\omega^n \cdot \displaystyle \int_M\;\omega^n\\
& = & \displaystyle \int_M\; (R - n )^2 \;\omega^n \cdot V
\\ &\leq & \epsilon_0^2 \cdot V \cdot V = \epsilon_0^2 \cdot V^2,
\end{array}
\]
which gives
\begin{equation}
{1 \over V} \displaystyle \int_M |{\rm Ric} - 1 |\;\omega^n \leq
\epsilon_0. \label{eq:dimaeter2}
\end{equation}
The value of $\epsilon_0$ will be determined later.

Using this inequality (\ref{eq:dimaeter2}), we want to show that
the diameter must be bounded from above. Note that in our
setting, $m = {\rm dim}(M) = 2 n.\;$ Unlike in \cite{spro001}, we
are not
interested in obtaining a sharp upper bound on the diameter. \\

Let $A_1$ and $A_2$ be two balla of small radius and $W = M. \;$
Let $f= |{\rm Ric} - 1| = \displaystyle \sum_{i=1}^m |\lambda_i -
1|,\;$ where $\lambda_i$ is the eigenvalue of the Ricci tensor.
We assume also that all geodesics are parameterized by arc length.
By possibly removing a set of measure $0$ in $A_1 \times   A_2,$
there is a unique minimal geodesic from $x$ to $y$ for all $(x,y)
\in A_1 \times A_2.$ Let $p,q$ be two points on $M$ such that
\[
   d(p,q) = {\rm diam}(M) = D.
\]

We also used $d\,vol$ to denote the volume element in the
Riemannian manifold $M$ and $V$ denote the total volume of $M.\;$
For $r> 0$, put $A_1 = B(p, r) $ and $A_2 = B(q,r).\;$ Then Lemma
3.1 implies that
\[
\begin{array} {l} \displaystyle \int_{A_1 \times A_2} \displaystyle
\int_{r_{x,y}} \; |{\rm Ric} - 1| \;d\,s\; d \; vol_{A_1 \times A_2} \nonumber \\
\qquad \qquad \leq C(n,k,R)(diam(A_2) vol(A_1) + diam(A_1)
vol(A_2)) \displaystyle \int_W \; |{\rm Ric} - 1| \;d\,vol.
\end{array}
\]
Taking infimum over both sides, we obtain
\begin{eqnarray} \displaystyle \inf_{(x,y) \in A_1 \times A_2 } \displaystyle
\int_{r_{x,y}} \; |{\rm Ric} - 1| \;d\,s\nonumber \\
\qquad \leq 2 \;r\; C(n)({1\over vol(A_1)} + {1\over  vol(A_2)})
\displaystyle \int_W \; |{\rm Ric} - 1|\;d\,vol
\nonumber \\
\leq 4 r C(n) {{D^n}\over r^n} {1\over V} \displaystyle \int_M \;
|{\rm Ric} - 1| \;d\,vol, \label{eq:CT2}
\end{eqnarray}
where the last inequality follows from the relative volume
comparison. We can then find a minimizing unit-speed geodesic
$\gamma$ from $x \in \overline{A_1}$ and $y \in \overline{A_2}$
which realizes the infimum, and will show that for $L =  d(x,y)$
much larger than $\pi,$ $\gamma$ can not be minimizing if the
right hand
side of (\ref{eq:CT2}) is small enough.\\

Let $E_1(t), E_2(t), \cdots E_m(t)$ be a parallel orthonormal
basis along the geodesic $\gamma$ such that $E_1(t) =
\gamma'(t).\;$ Set now $Y_i(t) = \sin \left({{\pi t}\over L}
\right) E_i(t), i =2,3,\cdots m.\;$ Denote by $L_i(s)$ the length
functional of a fixed endpoint variation of curves through
$\gamma$ with variational direction $Y_i,$ we have the 2nd
variation formula
\[
\begin{array}{lcl}
 & & \displaystyle
\sum_{i=2}^m {{d^2 L_i(s)}\over {d\,s^2}}\mid_{s=0} \\
& = & \displaystyle \sum_{i=2}^m \displaystyle \int_0^L \left(
g(\nabla_{\gamma'} Y_i, \nabla_{\gamma'}Y_i) -
R(\gamma',Y_i,\gamma',Y_i) \right)d\,s\\
& = & \displaystyle \int_0^L (m-1) \left({\pi^2 \over L^2} \cos^2
\left( {{\pi t}\over L} \right) \right)- \sin^2 \left( {{\pi
t}\over L}
\right) {\rm Ric}(\gamma',\gamma') \;d\,s \\
& = & \displaystyle \int_0^L \left((m-1) {\pi^2 \over L^2} \cos^2
\left( {{\pi t}\over L} \right) -  \sin^2 \left( {{\pi t}\over L}
\right) \right)\;d\,s \\
& & \qquad + \displaystyle \int_0^L\; \sin^2 \left( {{\pi t}\over
L} \right) \left( 1 - {\rm Ric}(\gamma',\gamma') \right)\;d\,s \\
& = & -{L\over 2} \left( 1- (m-1) {\pi^2 \over L^2} \right) \\
& & \qquad + \displaystyle \int_0^L\; \sin^2 \left( {{\pi t}\over
L} \right) \left( 1 - {\rm Ric}(\gamma',\gamma') \right)\;d\,s.
\end{array}
\]

Note that
\[
  1 -{\rm Ric}(\gamma',\gamma') \leq |{\rm Ric} - 1|.
\]

Combining the above calculation and the inequality
(\ref{eq:CT2}), we obtain

\begin{eqnarray} & & \displaystyle
\sum_{i=2}^n {{d^2 L_i(s)}\over {d\,s^2}}\mid_{s=0}  \nonumber \\
&\leq  & -{L\over 2} \left( 1- (m-1){\pi^2 \over L^2} \right)  +
\displaystyle \int_0^L\; \sin^2 \left( {{\pi t}\over L} \right)
|{\rm Ric} - 1|\;d\,s \nonumber\\
& \leq & -{L\over 2} \left(1- (m-1){\pi^2 \over L^2} \right)  + 4
r C(n) {{D^n}\over r^n} {1\over V} \displaystyle \int_M \; |{\rm
Ric} - 1| \;d\,vol. \label{eq:CT3}
\end{eqnarray}
Here in the last inequality, we have already used the fact that
$\gamma$ is a geodesic which realizes the infimum of the left
side of inequality (\ref{eq:CT2}).  For any fixed positive larger
number $N > 4$,  let $D= N\cdot r .\;$ Set $ c = {1 \over V}
\displaystyle \int_M \; |{\rm Ric} - 1| \;d\,vol.\;$  Note that
\[ L = d(x,y) \geq d(p,q) - 2 r = D (1- {2\over N})\geq {D\over 2}.
\]
Then the above inequality (\ref{eq:CT3}) leads to
\[
\begin{array}{lcl}
  { 1 \over D}  \displaystyle
\sum_{i=2}^n {{d^2 L_i(s)}\over {d\,s^2}}\mid_{s=0}
 & \leq &  -{{1 - {2 \over N} }\over 2} \left(1- (m-1) {\pi^2 \over
L^2} \right) + 4 C(n) \; {N^{m-1} \over V} \cdot c  \cdot V\\
& = &  4 C(n) \; {N^{m-1}} \left( c - {{(N-2)}\over {2N}}{1 \over
{4 C(n) N^{m-1}}} \right)  + {{1 - {2 \over N} }\over
2}(m-1){\pi^2 \over L^2}.
\end{array}
\]
Note that the second term in the right hand side can be ignored if
$L \geq {D \over 2} $  is large enough.
Set
\[
\epsilon_0  = {{(N-2)}\over {2N}} \cdot {1 \over {4 C(n) N^{m-1}}}
= {{N-2 }\over {8 C(n) N^m  }}.
\]
Then if
\[
 {1 \over V} \displaystyle \int_M \; (R -n)^2 \omega^n \leq
 \epsilon_0^2,
\]
by the argument at the beginning of this proof, we have the
inequality (\ref{eq:dimaeter2}):
\[
 {1 \over V} \displaystyle \int_M \; (R -n)^2 \omega^n \leq
 \epsilon_0^2,
\] \[
{1\over V} \displaystyle \int_M \; |{\rm Ric} - 1| \;d\,vol <
\epsilon_0,
\]
which in turns imply
\[
 { 1 \over D}  \displaystyle
\sum_{i=2}^n {{d^2 L_i(s)}\over {d\,s^2}}\mid_{s=0} < 0,
 \]
for $D$ large enough. Thus, if the diameter is too
large,$\;\gamma$ cannot be a length minimizing geodesic. This
contradicts our earlier assumption that $\gamma$ is a
 minimizing geodesic. Therefore, the diameter must have a uniform upper bound.
\end{proof}

According to the work of.
 C. Croke \cite{croke80}, Li-Yau \cite{liyau80} and Li
\cite{pli80}), we state the following lemma on the upper bound of
the Sobolev constant and Poincare constant:
\begin{lem}
Let $(M,\omega)$ be any compact polarised K\"ahler manifold where
$[\omega]$ is the canonical class. If $Ric(\omega) \geq 0,\; V =
\displaystyle \int_M \omega^n  \geq \nu
> 0$ and the diameter has a uniform upper bound,
then there exists a constant $\sigma = \sigma(\epsilon_0, \nu)$
such that for all function $f\in C^{\infty}(M),\;$ we have
\[
\left( \displaystyle \int_M\;\mid f\mid^{{2n}\over {n-1}}\;
\omega^n \right)^{{n-1}\over n} \leq \sigma \left( \displaystyle
\int_M \mid \nabla f\mid^2 \;\omega^n + \displaystyle \int_M f^2
\;\omega^n \right).
\]
Furthermore, there exists a uniform Poincar\'e constant
$c(\epsilon_0) $ such that the Poincar\'e inequality holds
\[
\displaystyle \int_M \left( f - {1\over V} \displaystyle \int_M
\;f \;\omega^n \right)^2 \;\omega^n \leq c(\epsilon_0)\;
\displaystyle \int_M \; \mid\nabla f\mid^2 \; \omega^n.
\]
Here $\epsilon_0$ is the constant appeared in Lemma 3.5.
\end{lem}
\begin{proof} 
Note that $(M,\omega)$ has a uniform upper bound on the diameter.
Moreover,  it has a lower volume bound and it has non-negative
Ricci curvature. Following a proof in \cite{pli80} which is based
on a result of C. Croke \cite{croke80}, we obtain a uniform upper
bound on the Sobolev constant (independent of metric!).

Recall a theorem of Li-Yau \cite{liyau80} which gives a positive
lower bound of the first eigenvalue in terms of the diameter when
Ricci curvature is nonnegative: \[ \lambda_1(\omega) \geq
{{\pi^2}\over {4 D^2}},
\]
here $\lambda_1,\; D$ denote the first eigenvalue and the diameter
of the K\"ahler metric $\omega.\;$ Now $D$ has a uniform upper
bound according to Lemma 3.5. Thus the first eigenvalue of
$\omega$ has a uniform positive lower bound; which, in turn,
implies that there exists a uniform Poincar\'e constant.
\end{proof}

\section{$C^0$ estimates }

Let us first prove a general lemma on $C^0$ estimate:
\begin{lem} Let $\omega_{\psi}$ be a K\"ahler metric such that
\[ \displaystyle \sup_M \psi \leq C_1,
\]
and
\[
  \displaystyle \int_M (- \psi) \omega_{\psi}^n \leq C_2.
\]
If the Sobolev constant and the Poincar\'e constant of
$\omega_{\psi}$ are bounded from above by $C_3,$ then there exists
a uniform constant $C_4$ which depends only on the dimension and
the constants $C_1, C_2 $ and $C_3\;$ such that
\[
     |\psi |\leq C_4.
\]
\end{lem}
We will use this lemma several times, so we include a proof here
for the convenience of the reader.
\begin{proof}
 Denote by $\Delta _{\psi}$ the Laplacian of $\o
_{\psi}.\;$ Then, because $\o + \p \b \p {\psi}
>0$, we see that $\o = \o _{\psi} - \p \b \p
{\psi} >0.\;$ Taking the trace of this latter expression with
respect to $\o _{\psi}$, we get
$$
n - \Delta _{\psi} {\psi} = \tr _{\o _{\psi}} \o
> 0.
$$

Define now ${\psi} _- (x) = \max \{ -{\psi} (x) , 1 \} \ge 1.\;$
It is clear  that
 \[
 {{\psi} _-} ^p (n - \Delta
_{\psi} {\psi} ) \ge 0. \] Integrating this inequality, we get
\begin{align*}
0 & \le {1 \over V} \int _M {{\psi} _- }^p (n - \Delta _{\psi} {\psi} ) \o _{\psi} ^n \\
  &= {n \over V} \int _M {{\psi} _- }^p \o _{\psi} ^n + {1 \over V} \int _M
     \nabla _{\psi} {{\psi} _- }^p \nabla _{\psi} {\psi} \o _{\psi} ^n \\
  &= {n \over V} \int _M {{\psi} _- }^p \o _{\psi} ^n + {1 \over V} \displaystyle \int _{\{ {\psi}
     \le -1 \} } \nabla _{\psi} {{\psi} _- }^p \nabla _{\psi} {\psi} \o _{\psi} ^n \\
  &= {n \over V} \int _M {{\psi} _- }^p \o _{\psi} ^n + {1 \over V} \int _M
      \nabla _{\psi} {{\psi} _- }^p \nabla _{\psi} (-{\psi} _- ) \o _{\psi} ^n \\
  &= {n \over V} \int _M {{\psi} _- }^p \o _{\psi} ^n  - {1 \over V} {4p \over
     {(p+1 )^2}} \int _M | \nabla _{\psi} {{\psi} _- }^{{{p+1} \over 2}} |^2
      \o _{\psi} ^n,
\end{align*}
which yields, using the fact that ${\psi} _- \ge 1$ and hence
${{\psi} _- }^p \le {{\psi} _- }^{p+1} ,$
$$
{1 \over V} \int _M \bigl| \nabla _{\psi} {{\psi} _-} ^{{{p+1}
\over 2}} \bigr|^2 \o _{\psi} ^n  \le {{n(p+1)^2 } \over {4p V}}
\int _M {{\psi} _- }^{p+1} \o _{\psi} ^n .
$$


Since the Sobolev constant of $\omega_{\psi}$ is bounded from
above, we can use the Sobolev inequality,
$$
{1 \over V} \Bl\int_M  |{\psi} _- | ^{{(p+1) n} \over {n-1}} \o
_{\psi} ^n \Br ^{{n-1} \over n} \le {{c(p+1)} \over V} \int _M
{{\psi} _- }^{p+1} \o _{\psi} ^n.
$$

Moser's iteration will show us that
$$
\sup _M {\psi} _- = \lim _{p \to \infty} \| {\psi} _- \| _{
L^{p+1} (M,\o _{\psi})} \le C \| {\psi} _- \| _{L^2 (M, \o
_{\psi} )}.
$$

Since the Poincare constant is uniformly bounded from above, we
can use the Poincar\'e inequality
\begin{align*}
{1 \over V} \int _M \Bl {\psi} _- - {C \over V} \int _M {\psi} _-
\o _{\psi} ^n \Br ^2 \o _{\psi} ^n
&\le {1 \over V} \int _M | \nabla {\psi} _- | ^2 \o _{\psi} ^n \\
&\le {C' \over V} \int _M {\psi} _- \o _{\psi} ^n,
\end{align*}
where we have set $p=1$ and used the same reasoning as before.
This then implies that
$$
\max  \{ - \inf _M {\psi} , 1 \} = \sup _M {\psi} _- \le {{C}
\over V} \int _M {\psi} _- \o _{\psi} ^n.
$$
Since $\int_M e^{-h_\varphi + {\psi}} \o_{\psi}^n = V$, we can
easily deduce $\displaystyle \int_{{\psi} > 0} {\psi} \o_{\psi}^n
\le C$. Combining this together with the above, we get
$$
-\inf _M {\psi} \le {C \over V} \int _M (-{\psi} ) \o _{\psi} ^n
+ C,
$$
which proves the lemma.
\end{proof}

\begin{lem} Along the K\"ahler Ricci flow, the diameter of the evolving metric is uniformly bounded.
\end{lem}
\begin{proof} In our first work \cite{chentian001}, we proved that
\[
 \displaystyle \int_0^\infty\;d\,t \displaystyle \int_M\; \; (R-r) ^2\;
 \omega_{\varphi}^{n}\leq C.
\]
Thefore, for any sequence $s_i \rightarrow \infty,$ and for any
fixed time period $T$, there exists $t_i \rightarrow \infty$ and
$0 < s_i -t_i < T$ such that
\begin{equation}
\displaystyle \lim_{t_i \rightarrow \infty} {{1}\over V}\;
\displaystyle \int_M\; (R-r) ^2\;   \omega_{\varphi}^{n}  =0.
\label{eq:sequence}
\end{equation}
Now for this sequence of $t_i, $ applying lemma 3.5, we show there
exists a uniform constant $D$ such that the diameters of
$\omega_{\varphi(t_i)}$ are uniformly bounded by ${D\over 2}.\;$
Recalled that the Ricci curvature is uniformly positive along the
flow so that diameter of evolving metric increased at most
exponentially since
\[
  {\p\over {\p t}} g_{i\bar j} = g_{i \bar j} - R_{i \bar j} \leq
  g_{i\bar j}.
\]
Now $t_{i+1} -t_i < 2 T$ for all $i>0,$ this implies that the
Diameters of the evolving metric along the entire flow is
controlled by $e^{2T} {D\over 2} \leq D$ (choose $T$ small enough
in the first place. \end{proof}
Combining this with Lemma 3.6, we obtain

\begin{theo} Along the K\"ahler Ricci flow,  the evolving K\"ahler metric $\omega_{\varphi(t)}$
has a uniform upper bound on the  Sobolev constant and Poincar\'e
constant.
\end{theo}

Before we go on any further,  we want to review some results we
obtained in previous paper \cite{chentian001}.

Let $\varphi(t)$ be the global solution of the K\"ahler Ricci
flow.  In the level of K\"ahler potentials, the evolution
equation is:
\[
{{\p \varphi(t)} \over {\p t}} = \log
\left({{\omega_{\varphi(t)}^n }\over {\omega^n}}\right) +
\varphi(t) - h_{\omega}.
\]

 According to Lemma 6.5 of \cite{chentian001}, there exists a one parameter
 family of
K\"ahler-Einstein metrics $\omega_{\rho(t)} = \omega + \sqrt{-1}
\partial \bar \partial \rho(t) $ such that $\omega_{\varphi(t)}$
is centrally positioned with respect to  $\omega_{\rho(t)} $ for
any $t \geq 0. \;$ Suppose that $\omega_{\varphi(0)}$ is already
centrally positioned with the K\"ahler-Einstein metric $\omega_1
= \omega + \sqrt{-1} \p \bar \p \rho(0).\;$ Normalize the value
of $\rho(t)$ such that
 $$\omega_{\rho(t)}^n = e^{-\rho(t) +
h_{\omega}} \; \omega^n,\;$$
 or equivalently
 \begin{equation} \ln \left(
{{\omega_{\varphi_{\rho(t)}}^n }\over {\omega^n}}\right) =
-\rho(t) + h_\omega. \label{eq:modifedke}
\end{equation}
 Then the
K\"ahler Ricci flow equation can be re-written as ,
\begin{equation}
 {{\p \varphi(t)} \over {\p t}}= \log
\left({{\omega_{\varphi(t)}^n }\over {\omega_{\rho(t)}^n}}\right)
+ \varphi(t) - \rho(t). \label{eq:modified}
\end{equation}
Sometimes we may refer this equation as the modified K\"ahler
Ricci flow.
Next we are ready to prove the $C^0$ estimates for both the
K\"ahler potentials and the volume form when $t=t_i.\;$

\begin{theo} 
There exists a uniform constant $C$ such that
\[
    \mid\varphi(t) - \rho(t) \mid < C,  \qquad {\rm and}\;\; |{{\p \varphi}\over
    {\p t}} | \leq C.
\]
In particular, we have
\[
    \mid \ln \det \left({{\omega_{\varphi(t)}^n} \over
    {\omega_{\rho(t)}^n}}\right) \mid < C.
\]
\end{theo}


We need a lemma on $L^1$  integral of the K\"ahler potentials.
\begin{lem} Along the K\"ahler Ricci flow on a K\"ahler Einstein
manifold,  there exists a uniform bound $C$ such that
\[
-  C \leq \displaystyle \int_M \left(\varphi(t) -\rho(t)\right)
\omega_{\varphi(t)}^n  \leq C.
\]
\end{lem}
\begin{proof}
As in section 11 of \cite{chentian001}, we define
\[
   c(t) = \displaystyle \int_M\; {{\p \varphi(t)}\over {\p t}} \;
   \omega_{\varphi(t)}^n.
\]
In a K\"ahler Einstein manifold, the K-energy has a uniform lower
bound along the K\"ahler Ricci flow. Thus
\[
 \displaystyle \int_{0}^\infty\; \displaystyle \int_M \; \mid\nabla  {{\p \varphi(t)}\over {\p t}}
  \mid_{\varphi(t)}^2 \; \omega_{\varphi(t)}^n \;d\,t \leq C.
\]
Therefore, we can normalize the initial value of K\"ahler
potential so that
\[
  c(0) = \displaystyle \int_{0}^\infty\; e^{-t}\; \displaystyle \int_M \; \mid\nabla  {{\p \varphi(t)}\over {\p t}}
  \mid_{\varphi(t)}^2 \; \omega_{\varphi(t)}^n \;d\,t \leq C.
\]
According to Lemma 11.1 of \cite{chentian001}, we have  $c(t)
> 0$ and
\[
\displaystyle \lim_{t\rightarrow \infty} c(t) = 0.
\]
In particular, this implies that there exists a constant $C$ such
that
\[\begin{array}{lcl}
  C &\geq & c(t) =\displaystyle \int_M\; {{\p \varphi(t)}\over {\p t}} \;
   \omega_{\varphi(t)}^n \\
    & = & \displaystyle \int_M\; \left( \log
\left({{\omega_{\varphi(t)}^n }\over {\omega_{\rho(t)}^n}}\right)
+ \varphi(t) - \rho(t)\right)\;
   \omega_{\varphi(t)}^n > 0.
\end{array}\]
In the last inequality we have used the fact that $c(t) > 0.\;$
According to Theorem 2.14, we have
\[
 - C \leq  \displaystyle \int_M \; \log \left( { {\omega_{\varphi}}^n \over {{\omega_{\rho(t)}}^n}}\right) \omega_{\varphi(t)}^n
  \leq  C.
\]
Combining this with the previous inequality, we arrive at
\[
- C \leq \displaystyle \int_M\; \left( \varphi(t) -
\rho(t)\right)\; \omega_{\varphi(t)}^n < C.
\]
Here $C$ is a constant which may be different from line to line.
\end{proof}

Next we return to the proof of Theorem 4.4.
\begin{proof} According to Proposition 2.14, we  have
\[ J_{\omega_{\rho(t)}}(\omega_{\varphi(t)})  < C.
\]
Then
\[
 0 \leq (I - J)(\omega_{\varphi(t)},\omega_{\rho(t)}) \leq (n+1)\cdot
 J_{\omega_{\rho(t)}}(\omega_{\varphi(t)}) < (n+1) C.
\]
By definition, this implies that
\[0 \leq
\displaystyle \int_M\; (\varphi(t) - \rho(t)) (\omega_{\rho(t)}^n
- \omega_{\varphi(t)}^n) \leq C.
\]
Combining this with Lemma 4.6 we obtain
\[
-C \leq \displaystyle \int_M\; (\varphi(t) - \rho(t))
\omega_{\rho(t)}^n
 \leq C.
\]
Since $\triangle_{\rho(t)} (\varphi(t) - \rho(t)) \geq -n,\;$ by
the Green formula, we have
  \[
  \begin{array}{lcl} & &
   \displaystyle \sup_M (\varphi(t)-\rho(t)) \\
   & \leq & {1\over V} \displaystyle \int_M (\varphi(t)- \rho(t)) \omega_{\rho(t)}^n
   - \displaystyle \max_{x \in M}\;  \left( {1\over V} \int_M (G(x,\cdot) + C_4)
   \triangle_{\rho(t)}(\varphi(t)-\rho(t))
   \omega_{\rho(t)}^n (y)\right)\\
   & \leq & {1\over V} \displaystyle \int_M (\varphi(t)-\rho(t)) \omega_{\rho(t)}^n + n
   C_4,
  \end{array}
  \]
  where $G(x,y)$ is the Green function associated to $\omega_{\rho}$ satisfying
  $G(x,\cdot)\geq 0.\;$ Therefore, there exists a uniform constant $C$ such that
  \[
  \displaystyle \sup_{M}\; (\varphi(t)-\rho(t)) \leq C.
  \]
  By Lemma 4.5, we have
  \[
-  C \leq \displaystyle \int_M \left(\varphi(t) -\rho(t)\right)
\omega_{\varphi(t)}^n  \leq C.
  \]
  Furthermore, according to Theorem 4.3, the K\"ahler metrics
  $\omega_{\varphi(t)}$ have a uniform upper bound on both the
  Sobolev constant and the Poincar\'e constant. Now using Lemma
  4.1, we conclude that there exists a uniform constant $C$ such
  that
\[
 - C \leq (\varphi(t)-\rho(t)) \leq C.
\]
 Next we consider the following
 \[
 \begin{array}{lcl}
 {{\p \varphi(t)}\over {\p t}}  & = & \log
\left({{\omega_{\varphi(t)}^n }\over
{\omega_{\rho(t)}^n}}\right) + (\varphi(t) - \rho(t))\\
& \geq & - C,
 \end{array}
 \]
 for some uniform constant $C.\;$  Recall that $ |c(t)| = |\int_M\;
 {{\p\varphi(t)}\over {\p t}} \;\omega_\varphi^n|\;$  is uniformly
 bounded.  Therefore, there is some uniform constant $C$ such
 that
 \[
   \int_M\; | {{\p \varphi}\over {\p t}}| \omega_\varphi^n \leq C.
 \]
 In view of the fact the K energy is
uniformly bounded below, we arrive at
\[ \int_0^\infty\;d\,t\int_M\; |\nabla {{\p \varphi}\over {\p t}}|^2 \;\omega_\varphi^n < \infty.
\]
Since  the Poincare constant of the evolving K\"ahler metric is
bounded,  we have
\[ \displaystyle\;\int_a^{a+1}\;d t\int_M\; \left({{\p \varphi}\over {\p t}}\right)^2 \;\omega_\varphi^n \leq
C,
 \]
where $c> 0$ is a constant independent of  $a> 0.\;$\\

 Note that ${{\p \varphi(t)}\over {\p t}}$ satisfies the following evolution equation \[
 {\p\over {\p t}} {{\p \varphi(t)}\over {\p t}}= \triangle_\varphi \; {{\p \varphi}\over {\p
 t}}  + {{\p \varphi}\over {\p
 t}}
 \]
 and the fact that both the Sobolev and the Poincare constants of the
 evolving metrics are uniformly bounded.
 Applying  Lemma 4.7, a parabolic version  of Lemma 4.1 below,
 we prove that there exists a uniform constant $C$ such that
 \[
   -C \leq {{\p \varphi(t)}\over {\p t}}  =  \log
\left({{\omega_{\varphi(t)}^n }\over {\omega_{\rho(t)}^n}}\right)
+ (\varphi(t) +  \rho(t)) < C.
 \]
It follows that
 \[
-C \leq  \log \left({{\omega_{\varphi(t)}^n }\over
{\omega_{\rho(t)}^n}}\right)  < C.
 \]
\end{proof}

By Proposition 2.14, there exists a one parameter family of
$\sigma(t) \in {\Aut}(M)$  such that $\omega_{\varphi(t)}$ is
centrally positioned with respect to the K\"ahler-Einstein metric
$\omega_{\rho(t)}.\;$ Here
\begin{equation}
\sigma(t)^* \omega_1 = \omega_{\rho(t)} = \omega + \sqrt{-1}
\partial \overline{\partial} \rho(t).
\label{eq:modifyautomorphism0}
\end{equation}

This condition ``centrally positioned" plays an important role in
deriving  Proposition 2.14 there. However, it is no longer needed
once we have Proposition 2.14.

\begin{lem} There exists a uniform constant $C$ such that
for all integers $i=1,2,\cdots, \infty, \;$ we have
\[
  |\rho(i) - \rho(i+1)|  < C.
\]
Moreover,
\[
   \mid \sigma(i+1) \sigma(i)^{-1} \mid_\hbar < C.
\]
Here $\hbar$ is the left invariant metric in ${\rm Aut}(M).\;$
\end{lem}
\begin{proof}
The modified K\"ahler Ricci flow is
\[
 {{\partial }\over {\partial t}} (\varphi - \rho) = \varphi - \rho
 + \log {{\omega_{{\varphi}}}^n \over {\omega_{\rho(t)}}^n } -  {{\partial  \rho}\over {\partial t}}.
\]
Since ${{\p \varphi}\over {\p t}}$ is uniformly bounded, we arrive
at
\[
\begin{array}{lcl} & &
  |\rho(i) - \rho(i+1)| \\
  & \leq & |\rho(i) -\varphi(i)| + |\rho(i+1) -
  \varphi(i+1) | + |\varphi(i) - \varphi(i+1)|\\
    & \leq & C.
    \end{array}
\]
Since $\omega_{\rho(t)}$ is a K\"ahler-Einstein metric for any
time $t,\;$ we have (cf. equation (\ref{eq:modifedke}))
\[
|\log {{\omega_{{\rho(i+1)}}}^n \over {\omega_{\rho(i)}}^n }| =
|\rho(i) - \rho(i+1)|<  C
\]
and
\[
 \mid \sigma(i+1) \sigma(i)^{-1} \mid_\hbar < C.
\]
\end{proof}

This lemma allows us to do the following modification on the curve
$\sigma(t) \in {\rm Aut} (M) $ (we also modify the curve $\rho(t)$
by the equation (\ref{eq:modifyautomorphism0})):  Fix all of the
integer points ($\sigma(i), i=1,2,\cdots $ )of the curve
$\sigma(t)$ first. At each unit interval, replace the original
curve in ${\rm Aut}(M)$ by a straight line which connects the two
end points in ${\rm Aut}(M).\;$  Such a new curve in ${\rm
Aut}(M)$ will satisfy all the estimates listed below (for
convenience, we still denote it as $\sigma(t), \rho(t)$
respectively):
\begin{enumerate}
\item Theorem 4.4 still holds for this new curve $\rho(t)$ since
we only change $\rho(t)$ by a uniformly controlled amount (fix at
each integer points, and adapt linear intepolation between them).
\item The new curves $\sigma(t)$ and $\rho(t)$ are Lipschitz with
a uniform Lipschitz constant for all the time $t \in [0, \infty).$
In fact, $\sigma(t)$ is a infinite long piecewise linear in ${\rm
Aut}_r(M).\;$ \item There exists a uniform constant $C$ such that
\[
 |\left({{d\,}\over {d\,t}} {\sigma}(t) \right) \cdot
  {\sigma}(t)^{-1} | < C,\qquad {\rm for\;\;any}\qquad t \neq \;\;
  {\rm integer}.
\]
\end{enumerate}

In the remaining of this section, we want to give a technical
lemma required by the proof of Theorem 4.4.

\begin{lem}\footnote{This is a parabolic verison of Moser iteration arguments. We give a detailed proof
here for the convenience of the readers.} If the Poincare constant
and the Sobolev constant of the evolving K\"ahler metrics are both
uniformly bounded along the K\"ahler Ricci flow, and if ${{\p
\varphi}\over {\p t}}$ is bounded from below  and if
$\displaystyle\;\int_a^{a+1}\;d t\int_M\; \left({{\p \varphi}\over
{\p t}}\right)^2 \omega_{\varphi(t)}^n $ is uniformly bounded from
above for any $a \geq 0,\;$ then ${{\p \varphi}\over {\p t}}$ is
uniformly bounded from above and below.
\end{lem}
\begin{proof}  Since ${{\p \varphi}\over {\p t}}$ has a uniform lower bound,
there is a constant $c$ such that $u = {{\p \varphi}\over {\p t}}
+ c
> 1$ holds all the time. Now $u$ satisfies the equation:

\[
{\p \over {\p t}}\; u = \triangle_\varphi \; u + u - c.
\]
Set $d\,\mu(t) = \omega_{\varphi(t)}^n$ as the evolving volume
element. Then
\[ {\p \over {\p t}}\; d\,\mu(t) = \triangle_\varphi u
\;d\,\mu(t).
 \]
  For any $a < b < \infty, $ define $\eta$ to
be any positive increasing function which vanishes at $a.\;$ Set
\[\psi(t,x) = \eta^2 u^{\beta-1}
\]
for any $\beta > 2.\;$  Then (Here $\p_t d\mu(t) =
\triangle_\varphi u d\,\mu(t)$)
\[\begin{array}{lcl} & &
  \int_a^b\;d\,t\int_M \; (\p_t u) \psi d\,\mu(t)\\
   & =  & \int_a^b\;d\,t \left(\p_t \int_M \; u \psi d\,\mu(t) - \int_M\; u\; {{\p \psi}\over {\p t}}
    - \int_M\; \psi \; u\; \triangle_\varphi\;u\right)\\
    & = & \eta(b)^2 \int_M \; u^{\beta} - \int_a^b\;d\,t \{\int_M
    \left( u 2 \eta \eta' u^{\beta-1} + u \eta^2 (\beta-1) u^{\beta-2} \p_t u \right)+\int_M\; \psi \; u\;
\triangle_\varphi\;u\}.
  \end{array}
\]
Thus,
\[
\begin{array}{lcl} && \int_a^b\;d\,t \left(\int_M \;\beta (\p_t u) \psi
d\,\mu(t) + \int_M\;\psi \; u\; \triangle_\varphi\;u\right) \\
& = &\eta(b)^2 \int_M \; u^{\beta} - \int_a^b\;d\,t \int_M
  2 \eta \eta' u^{\beta} \\
  & = & \int_a^b\;d\,t\left(\int_M\; \beta (\triangle_\varphi u + u -c) \eta^2 u^{\beta-1} + \int_M\;\eta^2 u^{\beta-1} \; u\;
  \triangle_\varphi\;u\right)\\
  & \leq & -\int_a^b\;d\,t\left(\int_M\; \beta (\beta-1)
  u^{\beta-2} |\nabla u|^2 \eta^2 - \int_M\;\beta \eta^2
  u^\beta\right).
\end{array}
\]
Therefore, we have
\[\eta(b)^2 \int_M \; u^{\beta} +
\int_a^b\;d\,t\int_M\; \beta (\beta-1)
  u^{\beta-2} |\nabla u|^2 \eta^2 \leq \int_a^b\;d\,t\int_M\;
  \beta \eta^2 u^\beta + \int_a^b\;d\,t \int_M
  2 \eta \eta' u^{\beta}.
\]
In other words,
\[
\begin{array}{lcl} &&
\eta(b)^2 \int_M \; u^{\beta} + \int_a^b\;d\,t\int_M\; 4 (1
-{1\over \beta})
   |\nabla u^{\beta \over 2}|^2 \eta^2 \\
  & \leq & \int_a^b\;d\,t\int_M\;
  \beta (\eta^2 + 2 \eta \eta') u^\beta
\end{array}
\]
or
\[
\begin{array}{lcl} &&
\eta(b)^2 \int_M \; u^{\beta} + \int_a^b\;d\,t\int_M\; 4 (1
-{1\over \beta}) \left(
   |\nabla u^{\beta \over 2}|^2 \eta^2 + u^\beta \eta^2\right) \\
  & \leq & \int_a^b\;d\,t\int_M\;
  \beta (2 \eta^2 + 2 \eta \eta') u^\beta.
\end{array}
\]
In particular, this implies that
\[
\max_{a\leq t \leq b} \int_M\;\eta(t)^2 u^\beta \leq
\int_a^b\;d\,t\int_M\;
  \beta (2 \eta^2 + 2 \eta \eta') u^\beta.
\]
Let us first state a lemma.
\begin{lem} (Sobolev inequality) Assume $0\leq a < b$ and $v: M
\times [a,b] \rightarrow {\bf R}$ is a measurable function such
that
\[
   \sup_{a\leq t \leq b}\; |v(\cdot, t)|_{L^2(M, d\mu(t))} < \infty
\]
and
\[
\int_a^b\; \int_M\; |\nabla v |^2 \;d\,\mu d\,t< \infty,
\]
then we have ($m=2n = \dim (M)$)
\[
\int_a^b\;d\,t\int_M |v|^{ 2(m+2) \over m }\;d\,\mu(t) \leq \sigma
\sup_{a\leq t \leq b}\; |v(\cdot, t)|_{L^2(M, d\mu(t))}^{4\over m}
\;\int_a^b \;d\,t\int_M\; \left(|\nabla v|^2 + v^2
\right)\;d\,\mu(t).
\]
Here $\sigma$ is the Sobolev constant.
\end{lem}
\begin{proof} For any $a\leq t\leq b$, we have
\[ \begin{array}{lcl} |v(\cdot,t)|_{L^{2(m+2) \over m}(M, d\mu(t))
}& \leq & |v(\cdot, t)|_{L^2(M,d\mu(t))}^{2\over {m+2}}
\;|v(\cdot,
t)|_{L^{2m\over {m-2}}(M,d\mu(t)) }^{m\over {m+2}}\\
& \leq & |v(\cdot, t)_{L^2(M,d\mu(t))}^{2\over {m+2}} \;
\left(\sigma\;\int_M\; \left(|\nabla v|^2 + v^2
\right)\;d\,\mu(t)\right)^{m\over {2(m+2)}}.
\end{array}
\]
The lemma follows by taking ${2(m+2)\over m}$ power on both sides
and integrating over $[a,b].\;$
\end{proof}

Now we return to the proof of main theorem. Let $v = \eta u^{\beta
\over 2}$, we have
\[\begin{array}{lcl} & & \left(\int_a^b\;d\,t \int_M\; |\eta^2 u^\beta|^{{m+2}\over m}\right)^{m\over
{m+2}}\\ & \leq &  \sigma^{m\over {m+2}} \sup_{a\leq t \leq b}\;
|v(\cdot, t)|_{L^2(M, d\mu(t))}^{4\over {m+2}} \;\left(\int_a^b
\;d\,t\int_M\; \left(|\nabla v|^2 + v^2
\right)\;d\,\mu(t)\right)^{m\over {m+2}} \\
& \leq & C(m) \left(\int_a^b\;d\,t\int_M\;
  \beta (2 \eta^2 + 2 \eta \eta') u^\beta \right)^{2\over {m+2}} \left(\int_a^b\;d\,t\int_M\;
  \beta (2 \eta^2 + 2 \eta \eta') u^\beta\right)^{m\over {m+2}}\\
& \leq & C(m) \beta \; \int_a^b\;d\,t\int_M\;
   ( \eta^2 +  \eta \eta') u^\beta.
\end{array}
\]
Here $C(m)$ is a constant depending only on the Sobolev constant
of $(M, g(t))$ and dimension of manifold.

Now for any $ a \leq b_0 < b \leq a+1,$ define
\[
  b_k = b -{{b-b_0}\over 2^k}
\]
for any $ k \in {\bf Z}_+.\;$  Fix a function $\eta_0 \in
C^\infty({\bf R}, {\bf R})$ such that $0 \leq \eta_0 \leq 1,
\eta'_0 \geq 0, \eta_0(t) = 0$ for $t \leq 0$ and $\eta_0(t) =1$
for $t\geq 1.\;$ For each integer $k > 0$, we let $\eta(t) =
\eta_0({{t-b_k}\over {b_{k+1}- b_k}})$ and $\beta = 2
\left({{m+2}\over m }\right)^k.\;$ Then ($b_{k+1} - b_k =
{{b-b_0}\over 2^{k+1}}$)
\[\begin{array}{lcl} &&
\left(\int_{b_{k+1}}^{a+1}\;d\,t \int_M\; u^{2 ({{m+2}\over
m})^{k+1}} d\,\mu \right)^{{1\over 2} ({{m}\over {m+2}})^{k+1}}
\\ & \leq & C(m)^{{1\over 2} ({{m}\over {m+2}})^{k}} \left({{m+2}\over
m}\right)^{k ({{m}\over {m+2}})^{k}}
\left({2^{k+1}\over{b-b_0}}\right)^{{1\over 2} ({{m}\over
{m+2}})^{k}} \left(\int_{b_k}^{a+1} \;d\,t \int_M\;u^{2
\left({{m+2}\over m}\right)^k} \;d\,\mu \right)^{{1\over 2}
\left({m\over {m+2}} \right)^k}.
\end{array}\]
The iteration shows that for any integer $k>0, $ we have
\[
\begin{array}{lcl} && \left(\int_{b}^{a+1}\;d\,t \int_M\; u^{2 ({{m+2}\over
m})^{k+1}} d\,\mu \right)^{{1\over 2} ({{m}\over {m+2}})^{k+1}}\\
&\leq & {C(m)\over (b-b_0)^{{m+2}\over 4}} \left(
\int_a^b\;d\,t\int_M\; u^2 d\,\mu \right)^{1\over 2}.
\end{array}
\]
Here again $C(m)$ is a uniform constant which depends only on the
sobleve constant of the evolving metrics and the dimension $m.\;$
Since the last term is uniformly bounded, this implies that as $k
\rightarrow \infty,$ we have
\[
\sup_{b\leq t\leq a+1} \; u \leq {C(m)\over (b-b_0)^{{m+2}\over
4}} \left( \int_a^b\;d\,t\int_M\; u^2 d\,\mu \right)^{1\over 2}.
\]
\end{proof}
\section{Uniform bounds  on gauge}

In order to use this uniform $C^0$ estimate and the flow equation
(\ref{eq:modified}) to derive  the desired $C^2$ estimate, we
still need to control the size of ${{\partial \rho} \over
{\partial t}}.\; $ However, from our earlier modification above,
we can not determine ${{\partial \rho} \over {\partial t}} $ at
any integer point. For any non-integer points, we have a uniform
bound $C$ such that
\[
\mid {{\partial \rho} \over {\partial
t}}\mid_t < C,\qquad \forall \; t \neq {\rm integer}.
\]

Note that $\sigma(t)$ is an infinite long broken line in ${\rm
Aut}_r(M). \;$ Next we want to further modify the curve
$\sigma(t)$ by smoothing the corner at the integer points. Let us
first set up some notations. Let
  $\jmath $ be the Lie algebra of ${\rm Aut}(M).\;$   As before, suppose $\hbar$ is the left invariant metric
on ${\rm Aut}(M).\;$
 Denote $id$ the identity element in ${\rm Aut}(M) $ and $\exp$ is the exponential
map at the identity.  Use $B_{r}$ to denote the ball centered at
the identity element with radius $r.\;$

After the modification of last section, $\sigma(t)$ is an
infinite long broken line in ${\rm Aut}(M).\;$ We can write down
this curve explicitly:  For any integer $i=0,1,2,\cdots, \infty,$
we have
 \begin{equation} \sigma(t) = \sigma(i) \cdot \exp((t-i)
 X_i),\qquad \forall\; t \in [i,i+1].
 \end{equation}
 Here $\{X_i\}$ is a sequence of vector fields in $\jmath$ with
 a uniform upper bound $C$ on their lengths:
 \begin{equation}
   \|X_i\|_{\hbar} \leq  C, \qquad \forall\; i=0,1,2,\cdots,
   \infty.
 \end{equation}
Then there exists a uniform positive number ${1\over 4} >\delta >
   0$ such that for any integer $i> 0,$ we have
   \[
     \sigma(t) \in \sigma(i) \cdot B_{1\over 2}, \qquad \forall \; t \in (i-\delta, i+\delta).
   \]

   Note that $\delta$ depends on $\|X_i\|_{\hbar}.\;$ Since the latter has a uniform upper bound, $\delta$ must have
   a uniform lower bound. We then can choose one $\delta > 0$ for
   all $i.\;$  \\

   At each ball $\sigma_i \cdot B_1,$ we want to replace the
   curve segment $\sigma(t) \;(t\in [i-\delta,i+\delta])$ by a new
   smooth curve $\tilde{\sigma}(t) $ such that:
   \begin{enumerate}
   \item The two end points and their derivatives are not changed\footnote{In a Euclidean ball, we can use the 4th order
   polynomial to achieve this. In any unit ball of any finite dimensional Riemannian manifold, we can always
   do this uniformly, as long as  the metric and other data involved are uniformly bounded. }:
   \[
     \tilde{\sigma}(i\pm\delta) = \sigma(i\pm\delta)  \]
     and
     \[
      \left(\left({{d\,}\over {d\,t}} \tilde{\sigma}(t)\right)  \tilde{\sigma}(t)^{-1}\right)_{t=i\pm\delta}
     =  \left(\left({{d\,}\over {d\,t}} \sigma(t)\right)  \sigma(t)^{-1}\right)_{t=i\pm\delta}.
     \]
   \item There exists a uniform bound $C'$ which depends only on
   the upper bound of $\|X_i\|_{\hbar}$ and $\delta$ such that
   \[
   \|\left({{d\,}\over {d\, t}} \tilde{\sigma}(t)\right)\tilde{\sigma}(t)^{-1}
   \|_{\hbar} \leq C', \qquad \forall\; t \in [i-\delta,i+\delta].
   \]
   \item For any $t\in [i-\delta, i+ \delta],\;$ we have $\tilde{\sigma} (t) \in \sigma(i)
   B_1.\;$In other words, there exists a uniform constant $C$ such
   that:
   \[
        |\tilde{\sigma}(t) {\sigma}(t)^{-1} |_{\hbar}  < C.
   \]
   \end{enumerate}
   The last step is to set $\tilde{\sigma}(t) = \sigma(t)$ for all
   other time. Then the new curve $\tilde{\sigma}(t)$ has all the
   properties we want:
\begin{enumerate}
\item There exists a uniform constant $C$ such that $|\tilde{\sigma}(t) {\sigma}(t)^{-1} | < C $
for all $ t \in [0,\infty).\;$
\item There  exists a uniform constant $C$ such that
\[
  |\left({{d\,}\over {d\,t}} \tilde{\sigma}(t) \right) \cdot
  \tilde{\sigma}(t)^{-1} | < C,\qquad {\rm for\;\;any}\qquad t \geq 0.
\]
\end{enumerate}

Denote by $\tilde{\sigma}(t)^* \omega_1 = \omega_{\tilde{\rho}(t)}
= \omega + \sqrt{-1} \partial \bar \partial \tilde{\rho}(t).\;$
Then, $\omega_{\tilde{\rho}(t)} $ is a K\"ahler-Einstein metric
\[
  \omega_{\tilde{\rho}(t)}^n = e^{-\tilde{\rho}(t) + h_{\omega}}
  \; \omega^n.
\]

There exists a uniform constant $C$ such that
\begin{equation}
|\rho(t) -\tilde{\rho}(t)| < C \label{eq:modified}
\end{equation}
and
\[
  |{{\partial \tilde{\rho}(t)} \over {\partial t}} | < C
\]
hold for all $t.\;$

Now the inequality (\ref{eq:modified}) implies that
\[
  |\log \left( \det {{\omega_{\rho(t)}^n} \over
  \omega_{\tilde{\rho}(t)}^n}\right)|  \leq C.
\]
Combining this with Theorem 4.4, we arrive at

\begin{theo} There exists a one parameter family of K\"ahler Einstein metrics
$\omega_{\tilde{\rho}(t)} = \omega + \sqrt{-1} \partial \bar
\partial \tilde{\rho}(t),$ which is essentially parallel to the initial family of K\"ahler Einstein metrics,
 and a uniform constant $C$ such that the following holds
\[
 |\varphi(t) - \tilde{\rho}(t)|  \leq  C,
\]

\[
 - C  < \log {{\omega_{{\varphi(t)}}}^n \over {\omega_{\tilde{\rho}(t)} }^n} <
 C
\]
and \[ |{{\partial \tilde{\rho}(t)} \over {\partial t}}| < C
\] over the entire modified K\"ahler Ricci flow.
\end{theo}

\section{$C^2$ and higher order derivative estimates}

Consider the modified K\"ahler Ricci flow
\begin{equation}
{{\partial }\over {\partial t}} (\varphi - \tilde{ \rho}) =
\varphi - \tilde{\rho}
 + \log {{\omega_{{\varphi}}}^n \over {\omega_{\tilde{\rho}(t)}}^n } -  {{\partial  \tilde{\rho}}\over {\partial t}}.
 \label{eq:newmodified}
\end{equation}
By Theorem 5.1, we have a uniform bound on both $\mid (\varphi -
\tilde{ \rho})\mid$ and $\mid{{\partial }\over {\partial t}}
(\varphi - \tilde{ \rho}) \mid .\; $ This fact will play an
important role in deriving $C^2$ estimate on the evolved relative
K\"ahler potential $(\varphi - \tilde{ \rho})$ in this section:

\begin{theo} If the $C^0$ norms of  $\mid
(\varphi - \tilde{ \rho})\mid$ and $\mid{{\partial }\over
{\partial t}} (\varphi - \tilde{ \rho}) \mid $ are uniformly
bounded (independent of time $t$), then  there exists a uniform
constant $C$ such that
\[
 0 \leq    n + \tilde{\triangle} (\varphi -\tilde{\rho}) < C,
\]
where $\tilde{\triangle}$ is the Laplacian operator corresponding
to the  evolved K\"ahler-Einstein metrics $\omega_{\rho(t)}. $
\end{theo}
 We set up some notations first. Let
$\triangle'$ be the Laplacian operator corresponding to the
evolved K\"ahler metric $\omega_{\varphi(t)}$ respectively.
 Let $\Box =  \triangle' - {{\partial} \over
{\partial t}}.\;$ Put $\omega_{\tilde{\rho}(t)} = \sqrt{-1}
h_{\alpha \bar \beta} d\,z^{\alpha} \otimes z^{\bar \beta}\;$ and
$ \omega_{\varphi(t)} = \sqrt{-1} g'_{\alpha \bar \beta}
d\,z^{\alpha} \otimes d\,z^{\bar \beta} \;$ where \[
g'_{\alpha\bar \beta} = h_{\alpha \bar \beta} + {{\p^2
\left(\varphi(t) -\tilde{\rho}(t) \right) }\over {\p z^{\alpha} \p
z^{\bar \beta}}}.
\] Then
\[
  \triangle' =  \displaystyle \sum_{\alpha,\beta=1}^n\; g'^{\alpha \bar \beta} {{\partial^2 }\over
  {\partial z^{\alpha} \partial z^{\bar \beta}}},\qquad
   \tilde{\triangle} = \displaystyle \sum_{\alpha,\beta=1}^n\; h^{\alpha \bar \beta} {{\partial^2\;}\over
  {\partial z^{\alpha} \partial z^{\bar \beta}}}.
\]
and
\[
[{{\partial} \over {\partial t}},  \tilde{\triangle}] = -
\displaystyle \sum_{a,b,c,d=1}^n\; h^{a \bar b} {{\partial^2\;
{{\partial \tilde{\rho}}\over {\partial t}}  }\over {\partial
z^{c} \partial z^{\bar b}}} h^{c \bar d} {{\partial^2 }\over
  {\partial z^{a} \partial z^{\bar d}}}.
\]
Furthermore, we have
\[
  \tilde{\triangle}  { {\partial \tilde{\rho}}\over {\partial
t}}  = - { {\partial \tilde{\rho}}\over {\partial t}}.
\]
Thus the  Hessian of ${ {\partial \tilde{\rho}}\over {\partial
t}}$ with respect to the evolved K\"ahler Einstein metric
$\omega_{\tilde{\rho}(t)}$ is uniformly bounded from above since
$|{ {\partial \tilde{\rho}}\over {\partial t}}|$ is uniformly
bounded from above.   \\

\begin{proof} of Theorem 6.1: We want to use the maximum principle in this proof.
Let us first calculate $ \Box \left( n + \tilde{\triangle}
(\varphi -\tilde{\rho}) \right).\;$

Let us choose a coordinate so that at a fixed point both
$\omega_{\tilde{\rho}(t)} =\sqrt{-1} h_{\alpha \bar \beta}
d\,z^{\alpha} \otimes d\,z^{\bar \beta}\; $  and the complex
Hessian of $\varphi(t)-\tilde{\rho}(t)$ are in diagonal forms. In
particular, we assume that $ h_{i \bar j} = \delta_{i \bar j} $
and $ \left(\varphi(t) -\tilde{\rho}(t)\right)_{i \bar j} =
\delta_{i \bar j} \left(\varphi(t) -\tilde{\rho}(t)\right)_{i \bar
i}.\;$ Thus
\[
g'^{ i \bar s} = {{\delta_{i \bar s}} \over {1 + \left(\varphi(t)
-\tilde{\rho}(t)\right)_{i \bar i} }}.\]

For convenience, put
\[
F = {{\partial }\over {\partial t}} \left(\varphi - \tilde{
\rho}\right) -\left(\varphi - \tilde{\rho}\right) + {{\partial
\tilde{\rho}}\over {\partial t}}.
\]

Note that $F$ has a uniform bound.  The modified K\"ahler Ricci flow (\ref{eq:newmodified})
can be reduced to
\[
 \log {{\omega_{{\varphi}}}^n \over {\omega_{\tilde{\rho}(t)}}^n }
 = F,
\]
 or, equivalently
 \[
 \left(\omega_{\tilde{\rho}(t)}  + \p \bar \p \left(\varphi- \tilde{\rho}\right) \right)^n
  = e^F \; \o ^n,
 \]
 i.e.,
 \[
 \log \det\left(h_{i \bar j} + {{\p^2 \left(\varphi - \tilde{\rho}\right)} \over {\p z_{i} \p z_{\bar j}}} \right)
 = F + \log \det (h_{i \bar j}).
 \]
 For convenience, set
 \[
    \psi(t) = \varphi(t) -\tilde{ \rho}(t)
 \]
 in this proof. Note that both $\mid \psi(t)\mid$ and $\mid {{\p \psi(t)} \over {\p t}}
 \mid$ are uniformly bounded (cf. Theorem 5.1).
 We first follow the standard calculation of $C^2$ estimates in \cite{Yau78}.
Differentiate both sides with respect to ${\p \over {\p z_k}}$
$$
(g')^{i \b j} \biggl( {{\p h_{i \b j}} \over {\p z_k}} +{{\p ^3
\psi(t)} \over { \p z_i \p \b z_j \p z_k}} \biggr) - h^{i \b j}
{{\p h_{i \b j}} \over {\p z_k}} = {{\p F} \over {\p z_k}},
$$
and differentiating again with respect to ${\p \over {\p \b
z_l}}$ yields
$$
(g')^{i \b j} \biggl( {{\p ^2 h_{i \b j}} \over {\p z_k \p \b
z_l}} + {{\p ^4 \psi(t)} \over { \p z_i \p \b z_j \p z_k \p \b
z_l}} \biggr) +h^{t \b j} h^{i \b s} {{\p h_{t \b s}} \over {\p
\b z_l}} {{\p h_{i \b j}} \over {\p z_k}} -h^{i \b j} {{\p ^2
h_{i \b j}} \over {\p z_k \p \b z_l}}
$$
$$
- (g')^{t \b j}(g')^{i \b s} \biggl( {{\p h_{t \b s}} \over {\p
\b z_l}} +{{\p ^3 \psi(t)} \over { \p z_t \p \b z_s \p \b z_l}}
\biggr)  \biggl( {{\p h_{i \b j}} \over {\p z_k}} +{{\p ^3
\psi(t)} \over { \p z_i \p \b z_j \p z_k}} \biggr)
 = {{\p ^2 F} \over {\p z_k \p \b z_l}}.
$$
Assume that we have normal coordinates at the given point, i.e.,
$h_{i \b j} = \delta _{ij}$ and the first order derivatives of
$g$ vanish. Now taking the trace of both sides results in
\begin{align*}
\Tilde{\Delta} F & = h^{k \b l} (g')^{i \b j}\biggl( {{\p ^2 h_{i
\b j}} \over {\p z_k \p \b z_l}} +{{\p ^4 \psi(t)} \over { \p z_i
\p \b z_j \p z_k \p \b z_l}}
 \biggr) \\ & \qquad -h^{k \b l}
(g')^{t \b j}(g')^{i \b s}{{\p ^3 \psi(t)} \over { \p z_t \p \b
z_s \p \b z_l}} {{\p ^3 \psi(t)} \over { \p z_i \p \b z_j \p
z_k}}  -h^{k \b l} h^{i \b j} {{\p ^2 h_{i \b j}} \over {\p z_k
\p \b z_l}} .
\end{align*}
On the other hand, we also have
\begin{align*}
\Delta ' (\Tilde{\Delta} \psi(t)) &= (g')^{k \b l} {{\p ^2 }
\over {\p z_k \p \b z_l}}
         \biggl( h^{i \b j} {{\p ^2 \psi(t)} \over {\p z_i \p \b z_j}}\biggr) \\
      &=(g')^{k \b l} h^{i \b j}  {{\p ^4 \psi(t)} \over { \p z_i \p \b z_j
          \p z_k \p \b z_l}} +(g')^{k \b l} {{\p ^2 h^{i \b j} } \over {\p z_k
         \p \b z_l}} {{\p ^2 \psi(t)} \over {\p z_i \p \b z_j}},
\end{align*}
and we will substitute ${{\p ^4 \psi(t)} \over {\p z_i \p \b
z_j\p z_k \p \b z_l}} $ in $\Delta ' (\Tilde{\Delta} \psi(t))$ so
that the above reads
\begin{align*}
\Delta ' (\Tilde{\Delta} \psi(t)) &= -h^{k \b l} (g')^{i \b j}
{{\p ^2 h_{i \b j}} \over
              {\p z_k \p \b z_l}} +h^{k \b l}
(g')^{t \b j}(g')^{i \b s}{{\p ^3 \psi(t)} \over { \p z_t \p \b
z_s \p \b z_l}}
{{\p ^3 \psi(t)} \over { \p z_i \p \b z_j \p z_k}} \\
 &\quad +h^{k \b l} h^{i \b j}
{{\p ^2 h_{i \b j}} \over {\p z_k \p \b z_l}} +\Tilde{\Delta} F
+(g')^{k \b l} {{\p ^2 h^{i \b j} } \over {\p z_k
         \p \b z_l}} {{\p ^2 \psi(t)} \over {\p z_i \p \b z_j}},
\end{align*}
which we can rewrite after substituting ${{\p ^2 h_{i \b j}} \over
{\p z_k \p \b z_l}} = -R_{i \b j k \b l}$ and  ${{\p ^2 h^{i \b
j}} \over {\p z_k \p \b z_l}} = R_{j \b i k \b l}$ as
$$
\begin{array}{lcl}
\Delta ' (\Tilde{\Delta} \psi(t)) & = & \Tilde{\Delta} F +h^{k \b
l}(g')^{t \b j}(g')^{i \b s} \psi(t) _{t \b s l} \psi(t) _{i \b j
k} \\ & & +(g')^{i \b j} h^{k \bar l} R_{i \b j k \b l} - h^{i \b
j} h^{k \bar l}R_{i \b j k \b l} +(g') ^{k \b l} R_{j \b i k \b l}
\psi(t) _{i \b j}.
\end{array}
$$
Restrict to  the coordinates we chose in the beginning so that
both $g$ and $\psi(t)$ are in diagonal form. The above transforms
to
$$
\begin{array}{lcl}
\Delta ' (\Tilde{\Delta} \psi(t)) & = & {1 \over {1 +\psi(t) _{i
\b i}}} {1 \over {1 +\psi(t) _{j \b j}}} \psi(t) _{i \b j k}
\psi(t) _{\b i j \b k} + \Tilde{\Delta} F \\ & & \qquad \qquad
+R_{i \b i k \b k} (-1+ {1 \over {1+\psi(t) _{i \b i}}}+{{\psi(t)
_{i \b i}} \over {1+\psi(t) _{k \b k}}}).
\end{array}
$$
Set now $C= \inf _{i \ne k} R_{i \b i k \b k}$ and observe that
\begin{align*}
 R_{i \b i k \b k} (-1+ {1 \over {1+\psi(t) _{i \b i}}}+{{\psi(t) _{i \b i}}
\over {1+\psi(t) _{k \b k}}})
  &= {1\over 2} { R_{i \b i k \b k}}{{(\psi(t) _{k \b k} -\psi(t) _{i \b i})^2}
     \over {(1 +\psi(t)_{i \b i})(1+\psi(t) _{k \b k})}} \\
  &\ge {C \over 2} {{(1+\psi(t) _{k \b k} -1-\psi(t) _{i \b i})^2} \over {(1 +\psi(t)
   _{i \b i})(1+\psi(t) _{k \b k})}} \\
  &= C\Bl {{1 +\psi(t) _{i \b i}} \over {1+\psi(t) _{k \b k}}} -1
  \Br,
\end{align*}
which yields
$$
\begin{array}{lcl}
\Delta ' (\Tilde{\Delta} \psi(t)) & \ge & {1 \over {(1 +\psi(t)
_{i \b i})(1 +\psi(t) _{j \b j})}} \psi(t) _{i \b j k} \psi(t)
_{\b i j \b k} + \Tilde{\Delta}
F \\
& & \qquad \qquad +C \Bl (n + \Tilde{\Delta} \psi(t)) \sum_i {1
\over {1 + \psi(t) _{i \b i}}}-1 \Br.
\end{array}
$$
We need to apply one more trick to obtain the requested
estimates. Namely,
\[
\begin{array} {lcl}
\Delta ' (e^{-\l \psi(t)} (n + \Tilde{\Delta} \psi(t)) )
  &= &  e^{-\l \psi(t)} \Delta ' (\Tilde{\Delta} \psi(t)) +2\nabla ' e^{-\l \psi(t)} \nabla '
   (n + \Tilde{\Delta} \psi(t)) \\
     & & \qquad +\Delta ' ( e^{-\l \psi(t)}) (n + \Tilde{\Delta} \psi(t)) \\
  &= & e^{-\l \psi(t)} \Delta ' (\Tilde{\Delta} \psi(t)) -\l  e^{-\l \psi(t)} (g')^{i \b i}
    \psi(t) _i (\Tilde{\Delta} \psi(t) )_{\b i}
     \\ & & \qquad -\l  e^{-\l \psi(t)} (g')^{i \b i} \psi(t) _{\b i}
     (\Tilde{\Delta} \psi(t)) _i \\
  & & \quad - \l  e^{-\l \psi(t)} \Delta ' \psi(t) (n + \Tilde{\Delta} \psi(t))
  \\ & & +\l ^2 e^{-\l \psi(t)}
    (g')^{i \b i} \psi(t) _{i}\psi(t)_{ \b i}  (n + \Tilde{\Delta} \psi(t)) \\
  &\ge &  e^{-\l \psi(t)} \Delta ' (\Tilde{\Delta} \psi(t))\\
  & &  -e^{-\l \psi(t)}(g')^{i \b i}
    (n + \Tilde{\Delta} \psi(t))^{-1}(\Tilde{\Delta} \psi(t) )_i (\Tilde{\Delta} \psi(t) ) _{\b i} \\
  & &\quad -\l e^{-\l \psi(t)}\Delta ' \psi(t) (n + \Tilde{\Delta} \psi(t)) ,
\end{array}
\]
which follows from the Schwarz Lemma applied to the middle two
terms. We will write out one term here, the other goes in an
analogous way
\begin{align*}
&(\l e^{-{\l \over 2} \psi(t)} \psi(t) _i (n+\Tilde{\Delta}
\psi(t) )^{{1 \over 2}} )
( e^{-{\l \over 2} \psi(t)} (\Tilde{\Delta} \psi(t)) _{\b i}(n+\Tilde{\Delta} \psi(t) )^{-{1 \over 2}})\\
\le & {1 \over 2} (\l ^2 e^{-\l \psi(t)} \psi(t) _i \psi(t) _{\b
i} (n+\Tilde{\Delta}
\psi(t) ) \\
& +e^{-\l \psi(t)}(\Tilde{\Delta} \psi(t)) _{\b i} (\Tilde{\Delta}
\psi(t))_i
   (n+\Tilde{\Delta} \psi(t) )^{-1}).
\end{align*}
Consider now the following
\begin{align*}
&-(n+\Tilde{\Delta} \psi(t) )^{-1} {1 \over {1+\psi(t) _{i \b
i}}} (\Tilde{\Delta} \psi(t)) _i (\Tilde{\Delta}
\psi(t) )_{\b i} +\Delta ' \Tilde{\Delta} \psi(t) \ge \\
& -(n+\Tilde{\Delta} \psi(t) )^{-1} {1 \over {1+\psi(t) _{i \b
i}}} |\psi(t) _{k \b k i}|^2
 +\Tilde{\Delta} F \\
&+ {1 \over {1+\psi(t) _{i \b i}}} {1 \over {1+\psi(t) _{k \b
k}}} \psi(t) _{k \b i \b j} \psi(t) _{i \b k j}
+C(n+\Tilde{\Delta} \psi(t) ){1 \over {1+\psi(t) _{i \b i}}}.
\end{align*}
On the other hand, using the Schwarz inequality, we have
\begin{align*}
& (n+\Tilde{\Delta} \psi(t) )^{-1}{1 \over {1+\psi(t) _{i \b i}}} |\psi(t) _{k \b k i}|^2 \\
& \quad \quad =(n+\Tilde{\Delta} \psi(t) )^{-1}{1 \over
{1+\psi(t) _{i \b i}}} \Biggl|{{\psi(t) _{k \b k i} } \over
{(1+\psi(t) _{k \b k} )^{1 \over 2}}}
 (1+\psi(t) _{k \b k} )^{1 \over 2} \Biggr| ^2 \\
&\quad \quad \le (n+\Tilde{\Delta} \psi(t) )^{-1} \Bl {1 \over
{1+\psi(t) _{i \b i}}}
 {1 \over {1+\psi(t) _{k \b k}}} \psi(t) _{k \b k i} \psi(t) _{\b k k \b i} \Br
\Bl 1+\psi(t) _{l \b l} \Br \\
&\quad \quad ={1 \over {1+\psi(t) _{i \b i}}}
 {1 \over {1+\psi(t) _{k \b k}}} \psi(t) _{k \b k i} \psi(t) _{\b k k \b i} \\
&\quad \quad ={1 \over {1+\psi(t) _{i \b i}}}{1 \over {1+\psi(t)
_{k \b k}}}
 \psi(t) _{i \b k k } \psi(t) _{k \b i \b k} \\
&\quad \quad \le {1 \over {1+\psi(t) _{i \b i}}}{1 \over
{1+\psi(t) _{k \b k}}}
 \psi(t) _{i \b k j } \psi(t) _{k \b i \b j},
\end{align*}
so that we get
$$
\begin{array}{l}
-(n+\Tilde{\Delta} \psi(t) )^{-1} {1 \over {1+\psi(t) _{i \b i}}}
(\Tilde{\Delta} \psi(t)) _i (\Tilde{\Delta} \psi(t) )_{\b i}
+\Delta ' \Tilde{\Delta} \psi(t) \\
\ge \Tilde{\Delta} F + C(n+\Tilde{\Delta} \psi(t) ) {1 \over {1 +
\psi(t) _{i \b i}}}.
\end{array}
$$
Putting all these together, we obtain

\begin{eqnarray}
 \triangle' \left(e^{-\lambda \psi(t) } (n + \tilde{\triangle} \psi(t)
 )\right) \qquad \qquad \nonumber \\
 \geq e^{-\lambda \psi(t) }
 \left( \tilde{\triangle} F + C  (n + \tilde{\triangle} \psi(t))
  \displaystyle \sum_{i=1}^n\;{1 \over {1 + \psi(t)_{i \bar i} }}
  \right) \nonumber \\
   \qquad \qquad - \lambda\; e^{-\lambda \psi(t)
  } \triangle' \psi(t) \; (n + \tilde{\triangle} \left(\varphi -
  \tilde{\rho}\right)). \qquad
  \label{eq:estimate1}
\end{eqnarray}

Consider
\[
\begin{array}{lcl} \tilde{\triangle} F  & = & \tilde{\triangle} \left({{\partial }\over {\partial t}} \psi(t) -\psi(t) + {{\partial
\tilde{\rho}}\over {\partial t}} \right) \\
& = & \tilde{\triangle} {{\partial }\over {\partial t}} \psi(t)
-( n+ \tilde{\triangle} \psi(t)) + n + \tilde{\triangle}
{{\partial \tilde{\rho}}\over {\partial t}} \\
& \geq & {{\partial }\over {\partial t}} ( n + \tilde{\triangle}
\psi(t) ) + \displaystyle \sum_{a,b,c,d=1}^n\; h^{a \bar b}
{{\partial^2\; {{\partial \tilde{\rho}}\over {\partial t}} }\over
{\partial z^{c} \partial z^{\bar b}}} h^{c \bar d} {{\partial^2
\psi(t) }\over
  {\partial z^{a} \partial z^{\bar d}}} \\
  & & \qquad \qquad \qquad \qquad \qquad \qquad- ( n+ \tilde{\triangle}
\psi(t)) + n + \tilde{\triangle}
{{\partial \tilde{\rho}}\over {\partial t}} \\
& \geq &  {{\partial }\over {\partial t}} ( n + \tilde{\triangle}
\psi(t) ) - c_1 ( n + \tilde{\triangle} \psi(t) ) - c_2
\end{array}
\]
for some uniform constants $c_1 $ and $c_2.\;$ In the last
inequality, we have used the fact that $|{{\partial
\tilde{\rho}}\over {\partial t}}|$ is uniformly bounded and
\[
\mid {{\partial^2\; {{\partial \tilde{\rho}}\over {\partial t}} }\over
{\partial z^{c} \partial z^{\bar b}}} \mid_{\tilde{\rho}(t)}   \leq c_3 \cdot
\mid {{\partial \tilde{\rho}}\over {\partial t}}\mid,
\]
and
\[
0 < h_{c\bar d} + {{\partial^2 \psi(t)}\over
  {\partial z^{c} \partial z^{\bar d}}} \leq \left(n +\tilde{\triangle} \psi(t) \right)
  \;\;h_{c\bar d}
\]
holds as matrix.  Here $c_3$ is some uniform constant.
\[
\begin{array}{lcl} e^{-\lambda \psi(t) } \tilde{\triangle} F & \geq
& e^{-\lambda \psi(t) } {{\partial }\over {\partial t}} ( n +
\tilde{\triangle} \psi(t) ) \\ & & \qquad \qquad - c_1
\;e^{-\lambda \psi(t) } \; ( n + \tilde{\triangle} \psi(t) ) -
c_2\; e^{-\lambda
\psi(t) } \\
& \geq &  {{\partial }\over {\partial t}}  \left( e^{-\lambda
\psi(t) }  ( n + \tilde{\triangle} \psi(t) ) \right) + \lambda
{{\partial }\over {\partial t}} \psi(t) e^{-\lambda \psi(t) }  ( n
+ \tilde{\triangle} \psi(t) )\\ & & \qquad \qquad - c_1
\;e^{-\lambda \psi(t) } \; ( n + \tilde{\triangle} \psi(t) ) -
c_2\; e^{-\lambda
\psi(t) } \\
& \geq & {{\partial }\over {\partial t}}  \left( e^{-\lambda
\psi(t) }  ( n + \tilde{\triangle} \psi(t) ) \right) -(c_1 +
|\lambda| c_4)\;e^{-\lambda \psi(t) } \; ( n +
\tilde{\triangle} \psi(t) ) \\
 & & \qquad \qquad \qquad- c_2\;
e^{-\lambda \psi(t) } .
\end{array}
\]
Here $c_4$ is a uniform constant.  Plugging this into the
inequality (\ref{eq:estimate1}), we obtain
\[ \begin{array}
{l}\Box   \left( e^{-\lambda \psi(t) } ( n + \tilde{\triangle}
\psi(t) )
\right) \\
\geq e^{-\lambda \psi(t) }
 \left(  C  (n + \tilde{\triangle} \psi(t))
  \displaystyle \sum_{i=1}^n\;{1 \over {1 + \psi(t)_{i \bar i} }}
  \right) \nonumber \\
   \qquad \qquad - \lambda\; e^{-\lambda \psi(t)
  } \triangle' \psi(t) \; (n + \tilde{\triangle} \psi(t)) \\  \qquad \qquad -(c_1 + |\lambda|
c_4)\;e^{-\lambda \psi(t) } \; ( n + \tilde{\triangle} \psi(t) )
- c_2\; e^{-\lambda \psi(t) }.
\end{array}\]

Now
\[ \triangle' \psi(t) =
n  - \displaystyle \sum_{i=1}^n\;{1 \over {1 + \left(\varphi -
\tilde{\rho}\right)_{i \bar i} }}.\]

 Plugging this into the
above inequality, we obtain
\[\begin{array}{l} \Box   \left( e^{-\lambda \psi(t) }
( n + \tilde{\triangle} \psi(t) ) \right) \\ \geq e^{-\lambda
\psi(t) }
 \left(  (C + \lambda)  (n + \tilde{\triangle} \psi(t))
  \displaystyle \sum_{i=1}^n\;{1 \over {1 + \psi(t)_{i \bar i} }}
  \right) \nonumber \\\qquad -(c_1 + |\lambda|
c_4 + n)\;e^{-\lambda \psi(t) } \; ( n + \tilde{\triangle}
\psi(t) ) - c_2\; e^{-\lambda \psi(t) }.
\end{array}
\]
Let $\lambda = - C + 1,$ we then have
 \[ \begin{array}{l}\Box   \left( e^{-\lambda \psi(t) }
( n + \tilde{\triangle} \psi(t) ) \right) \\ \geq e^{-\lambda
\psi(t) }
 \left( (n + \tilde{\triangle} \psi(t))
  \displaystyle \sum_{i=1}^n\;{1 \over {1 + \psi(t)_{i \bar i} }}
  \right) \nonumber \\\qquad  \qquad - c_5\;e^{-\lambda \psi(t) } \; ( n +
\tilde{\triangle} \psi(t) ) - c_2\; e^{-\lambda \left(\varphi -
\tilde{\rho}\right) }.
\end{array}\]
Here $c_5$ is a uniform constant.
 Note the following algebraic inequality
$$
\begin{array}{lcl}
\sum _i {1 \over {1 + \psi(t) _{i \b i}}} &  \ge & \biggr( {{\sum
_i ( 1 + \psi(t) _{i \b i})} \over {\prod _i(1 + \left(\varphi
- \tilde{\rho}\right) _{i \b i})}}\biggr) ^{1 \over {n-1}} \\
& = & e^{-{F \over {n-1}}} (n+\tilde{\Delta}\psi(t) ) ^{1 \over
{n-1}}.
\end{array}
$$

This can be verified by taking the $(n-1)$-th power of both
sides. So the last term in the above can be estimated by
$$\begin{array}{l}
e^{-\l \psi(t)}\sum _i {1 \over {1 + \psi(t)_{i \b i}}}(n + \Delta \psi(t))\\
\qquad \qquad
 \ge e^{-{F \over {n-1}}}
e^{-{\l \over {n-1}}} (e^{-\l \psi(t)} (n + \tilde{\Delta}
\psi(t))) ^{n \over {n-1}} .
\end{array}
$$
Setting now $u=e^{-\l \psi(t)} (n + \tilde{\Delta} \psi(t))$ and
recalling that $\psi(t) \le -1$ and hence $e^{-\l \psi(t)} \ge
1$,we finally obtain the following estimate
$$
\Box u \ge -c_1 -c_2 u +c_0 u^{n \over {n-1}}.
$$
Assume that $u$ achieves its maximum at $x_0$ and ${{\partial u}
\over {\partial t}} \mid_{x_0,t} \geq   0$, then at this point,
$\Box u =  \Delta ' u  - {{\partial u} \over {\partial t}}
\mid_{x_0,t}\le 0\;$ and therefore the maximum principle gives us
an upper bound $u(x_0 ) \le C$ which, in  turn, gives
 \[
0 \le \left(n + \tilde{\Delta} \psi(t) \right) (x) \le e^{\l
\psi(t) (x)} u(x_0 ) \le C
 \]
 and hence we found a $C^2$-estimate of $\psi(t).$

\end{proof}
\begin{prop}
Let $\tilde{\rho}(t) $ be as in Theorem 6.1.  Then there exists a
uniform constant $C$ such that
\[
\|\varphi(t) -\tilde{\rho}(t)\|_{C^3(\omega_{\tilde{\rho}})} \leq
C.
\]

\end{prop}
\begin{proof}
  Let
\[
   g'_{i \bar j} = h_{i \bar j} + \left(\varphi - \tilde{\rho}\right)_{i \bar j}
\]
and
\[
S = \displaystyle \sum_{i,j,k,r,s,t=1}^n\; g'^{i \bar r} g'^{\bar
j \;s} g'^{k\bar t} \left(\varphi - \tilde{\rho}\right)_{i \bar j
\;k} \left(\varphi - \tilde{\rho}\right)_{\bar r\; s \bar t}.
\]

 Using Calabi's computation and Theorem 5.1 as in \cite{Yau78}, one
 can show that $S\leq C$ for some uniform constant $C$.
 Consequently, the proposition is proved.
 \end{proof}

 \section{The proof of main theorems}
According to Theorems 5.1, 6.1 and 6.2, we have uniform $C^3$
estimates on $\varphi(t) - \tilde{\rho}(t)$  along the modified
K\"ahler Ricci flow. It is not difficult to prove the following
\begin{lem} For any integer $l>0$, there exists a uniform constant
$C_l$ such that
\[
   \| D^l\; \left(\varphi(t) - \tilde{\rho}(t)
   \right)\|_{\omega_{\tilde{\rho}}} \leq C_l,
\]
where $D^l$ represents arbitrary $l-$th derivatives.
Consequently,  there exists a uniform bound on the sectional
curvature and all the derivatives of $\omega_{\varphi(t)}.\;$ The
bound may possibly depend on the order of derivatives.
\end{lem}

Follow this lemma, we can easily derive that the evolved K\"ahler
metrics $\omega_{\varphi(t)} $ converge to a K\"ahler metric in
the limit (by choosing subsequence). We would like to show that
the limit is a K\"ahler-Einstein metric.  Following proposition
2.5 and the fact that $E_0$ and $E_1$ have a uniform lower bound,
we have
\[
\begin{array}{lcl}
 \displaystyle \int_0^\infty\; {{n\sqrt{-1}}\over V} \displaystyle \int_M \partial {{\partial
\varphi}\over {\partial t}} \wedge \overline{\partial} {{\partial
\varphi}\over{\partial
t}}  {\omega_{\varphi}}^{n-1}\; d\,t & = & E(0) - E(\infty) < C,\\
  \displaystyle
\int_0^\infty\;{{2}\over V} \displaystyle \int_M
(R(\omega_\varphi)-r)^2 {\omega_{\varphi}}^{n}\, d\,t & = &
E_1(0)-E_1(\infty)  \leq C. \end{array}
\]

Combining this with Lemma 7.1, we prove that for almost all
convergence subsequence of the evolved K\"ahler metrics
$\omega_{\varphi(t)},\;$ the limit metric is of constant scalar
curvature metric. From here, it  is not difficult to show that any
sequence of the evolved  K\"ahler metrics will have a subsequence
which converges to a metric of constant scalar curvature. In the
canonical class, any metric of constant scalar curvature is a
K\"ahler Einstein metric.
We then prove the following
\begin{theo} The modified K\"ahler Ricci flow converges to some
K\"ahler Einstein metric  by taking sub-sequences.
\end{theo}

To prove uniqueness of the limit by sequence, we can follow
\cite{chentian001} to first prove the exponential decay of
\[
  \displaystyle \int_M \left( {{\p \varphi}\over {\p t}} \right)^2
  \;\omega_{\varphi}^n.
\]
In other words, there exists a positive constant $\alpha$ and a
uniform constant $C$ such that
\[
 \displaystyle \int_M \left( {{\p \varphi}\over {\p t}} \right)^2
  \;\omega_{\varphi}^n < C e^{-\alpha t}
\]
for all evolved metrics over the K\"ahler Ricci flow. Eventually,
we prove the following main proposition (like in
\cite{chentian001}):
  \begin{prop} For any integer $l>0,\;{{\partial \varphi} \over {\partial t}} $
   converges exponentially fast to $0\;$ in any $C^l$ norm. Furthermore, the
   K\"ahler Ricci flow converges exponentially fast to a unique
   K\"ahler Einstein metric on any K\"ahler-Einstein manifolds.
   \end{prop}
\section{K\"ahler-Einstein orbifolds}
In this section, we will prove that any K\"ahler-Einstein orbifold
such that there is another  K\"ahler metric in the same K\"ahler
class  which has strictly positive bisectional curvature must be a
global quotient of $\CC P^n $ by a finite group. The simplest
example of K\"ahler orbifolds is the global quotient of $\CC P^n $
by a finite group. Roughly speaking, a generic K\"ahler orbifold
is the union of a family of open sets, where each open set admits
a finite  covering from an open smooth K\"ahler manifold where a
finite group acts holomorphically (we will give precise definition
later). If it admits a K\"ahler Einstein metric, then it is called
a K\"ahler-Einstein orbifold. The goal in this section is to show
that under our assumption, there exists a global branching
covering with a finite group action from $\CC P^n$ to the
underlying K\"ahler orbifold. The organization of this section is
as follows: In  subsection 8.1, we introduce the notion of complex
orbifolds and various geometric structures associated with them.
In Subsection 8.2, we consider the K\"ahler Ricci flow on any
K\"ahler Einstein orbifolds. If there is another K\"ahler metric
in the same K\"ahler class such that the bisectional curvature is
positive, then the K\"ahler Ricci flow converges and the limit
metric is a K\"ahler-Einstein metric with positive bisectional
curvature. In Subsection 8.3, we prove that any orbifold which
admits a K\"ahler-Einstein metric of constant bisectional
curvature  must be a global quotient of $\CC P^n.\;$ In subsection
8.4, we re-prove that any K\"ahler-Einstein metric with positive
bisectional curvature must be of constant bisectional curvature
(Berger's theorem \cite{Berger65}). We also prove that if a
K\"ahler metric is sufficiently close to a K\"ahler Einstein
metric on the K\"ahler Ricci flow, then the positivity of
bisectional curvature will be preserved when taking limit (Lemma
8.20).

\subsection{K\"ahler orbifolds}

Let us begin with the definition of  uniformization system over
an open connected analytic space \footnote{One reference for
orbifolds is Ruan \cite{Yuan00}. }:

\begin{defi} Let $U$ be a connected analytic space
 and $V$  a connected $n-$dimensional smooth K\"ahler manifold and
$G$  a finite group acting on $V\;$ holomorphically. An
$n-$dimensional uniformization system of $U$ is a triple $(V, G,
\pi)$,  where $\pi:\; V \rightarrow U$ is an analytic map inducing
an identification  between two analytic spaces $V/G$ and $U.\;$
Two uniformization systems $(V_i, G_i, \pi_i),\; i=1, 2, $ are
isomorphic if there is a bi-holomorphic map $\phi: V_1
\rightarrow V_2$ and isomorphism $\lambda: G_1 \rightarrow G_2$
such that $\phi$ is $\lambda-$equivariant, and $\pi_2 \circ \phi =
\pi_1.\;$
\end{defi}
In the above definition, we require that the fixed point set to be
real codimension 2 or higher (if the group action preserves
orientation, then the fixed point must be codimension 2 or
higher.). Then the non-fixed point set (the complement of the
fixed point set) is locally connected, which  is important for our
purpose. The following proposition is immediate:
\begin{prop} Let $(V, G, \pi)$ be a uniformization system of
$U.\;$For any connected open subset $U'$ of $U, \;(V, G, \pi)$
induces a unique isomorphism class of uniformization systems of
$U'.\;$
\end{prop}
\begin{proof} We want to clarify what ``induces" means in this
proposition. For any open subset $U' \subset U, $ consider the
preimage $\pi^{-1}(U')$ in $V.\;$  $G$ acts as permutations on the
set of connected components of $\pi^{-1}(U').\;$ Let $V'$ be one
of the connected components of $\pi^{-1}(U'),\; G' $ the subgroup
of $G$ which fixes the component $V'$ and $\pi'= \pi\mid_{V'}.\;$
Then $(V', G', \pi')$ is an induced uniformizing system of $U'.\;$
One can also show that any other induced uniformization system
must be isomorphic to this one.  We skip this part of the proof
and refer interested readers to \cite{Yuan00} for details.
\end{proof}

In light of this proposition, we can define equivalence of two
uniformization systems at a single point: For any point $p \in U,$
let $(V_1, G_1, \pi_1)$ and $(V_2, G_2, \pi_2)$ be two
uniformization systems of neighborhoods $U_1$ and $U_2$ of $p.\;$
We say that $(V_1, G_1, \pi_1)$ and $(V_2, G_2, \pi_2)$ are {\it
equivalent} at $p$ if they induce isomorphic uniformization
systems for a smaller neighborhood $U_3 \subset U_1 \bigcap U_2 $
of $p.\;$ Next we define a complex (K\"ahler) orbifold.

\begin{defi} Let $M$ be a connected analytic space.
 An $n-$dimensional complex orbifold structure on $M$ is
given by the following data: for any point $p \in M, $ there are
neighborhoods $U_p$ and their $n-$dimensional uniformization
systems $(V_p, G_p, \pi_p)$  such that for any $q \in U_p, $
$(V_p, G_p, \pi_p)$ and $(V_q, G_q,\pi_q)$ are equivalent at
$q.\;$  A point $p \in M$ is called regular if there exists a
uniformization system $(V_p, G_p, \pi_p)$ over $U_p \ni p $ such
that $G_p$ is trivial; Otherwise it is called singular. The set of
regular points is denoted by $M_{reg}.\;$ The set of singular
points is denoted by $M_{sing},\;$ and $M = M_{reg} \bigcup
M_{sing}.\;$
\end{defi}

 Next we define orbifold vector bundles over a complex orbifold. As before, we begin
 with local uniformization systems for orbifold vector bundles. Given
 an analytic space $U$ which is uniformized by $(V,G,\pi)$  and a complex analytic space $E$
 with a surjective holomorphic  map $pr: E \rightarrow U, $ a
 uniformization system of rank $k$ complex vector bundle for $E$ over $U$
 consists of the following data.
 \begin{enumerate}
 \item A uniformization system $(V, G, \pi)$ of $U.\;$
 \item A unifomization system $(V \times \CC^k, G, \tilde{\pi})$
 for $E.\;$ The action of $G$ on $ V \times \CC^k$ is an extension
 of the action of $G$ on $V$ given by $g(x, v) =(g \cdot x,
 \rho(x,g) \cdot v),$ where $\rho: V \times G \rightarrow
 GL(\CC^k)$ is a holomorphic map satisfying
 \[
  \rho(g \cdot x, h) \circ \rho(x, g) = \rho(x, h\circ g), \qquad
  \forall g,\;h \in G,\; x \in V.
 \]
 \item The natural projection map $\tilde{pr}: V \times \CC^k
 \rightarrow V$ satisfying
 \[
   \pi \circ \tilde{pr} = pr \circ \tilde{\pi}.
 \]
 \end{enumerate}
We can similarly define isomorphisms between two uniformization
systems of orbifold vector bundles for $E$ over $U.\;$ The only
additional requirement is that the diffeomorphism between $V
\times \CC^k$ are  linear on each fiber of $ \tilde {pr}: V \times
\CC^k \rightarrow V.\;$  Moreover, we can also define the
equivalent relation between two uniformization systems of complex
vector bundles at any specific point. Here is the definition of
orbifold vector bundles over complex orbifolds:

\begin{defi} Let $M$  be a complex orbifold and $E$ a complex vector space
with a surjective holomorphic  map $pr: E \rightarrow M.\;$ A rank
$k$ complex orbifold vector bundle structure on $E$ over $M$
consists of the following data: for each point $p\in M, $ there
is a unformized neighborhood $U_p$ and a uniformization system of
a rank $k$ complex vector bundle for $pr^{-1}(U_p)$ over $U_p$
such that for any $q \in U_p,$ the rank $k$ complex orbifold
vector bundles over $U_p$ and $U_q$ are isomorphic in a smaller
open subset $U_p \bigcap U_q.\;$ Two orbifold vector bundles
$pr_1: E_1 \rightarrow M$ and $pr_2: E_2 \rightarrow M$ are
isomorphic if there is a holomorphic map $\tilde{\psi}: E_1
\rightarrow E_2$ given by $\tilde{\psi}_p: (V_{1,p} \times \CC^k,
G_{1,p}, \tilde{\pi}_{1,p}) \rightarrow (V_{2,p} \times \CC^k,
G_{2,p}, \tilde{\pi}_{2,p}) $ which induces an isomorphism between
$(V_{1,p}, G_{1,p}, \tilde{\pi}_{1,p}) $ and $(V_{2,p}, G_{2,p},
\tilde{\pi}_{2,p}),$ and is a linear isomorphism between the
fibers of $\tilde{pr}_{1,p}$ and $\tilde{pr}_{2,p}.\;$
\end{defi}
For a complex orbifold, one can define the tangent bundle, the
cotangent bundle, and various exterior or tensor powers of these
bundles. All the differential geometric quantities such as
cohomology class, connections, metrics, and curvatures can be
introduced on the complex orbifold. \\

Suppose $M$ is a complex orbifold as in Definition 8.3. For any
$p\in M, $ let $p \in U_p$ be uniformized by $(V_p,G_p,\pi_p).\;$
When we say a metric $g$ is defined on $U_p, $ we really mean a
metric $\overline{g}$  defined on $V_p$ such that $G_p$ acts on
$V_p$ by isometries. For simplicity, we say the metric $g$  is
defined on $U_p$ and $\pi_p^* g = \overline{g}.\;$ This
simplification makes sense especially when $p$ is a regular point,
i.e., when $ G_p$ is trivial. One way to define a metric on the
entire complex orbifold is first to define it on $M_{reg}, $ then
extend it to be a metric on $M$ with possible singularities  since
$M_{sing}$ is codimension at least 2 or higher. The following
gives a definition of what a smooth K\"ahler metric or a K\"ahler
form on the complex orbifold is:

\begin{defi} For any point $p \in M, $ let $U_p$ be uniformized by $(V_p, G_p,
\pi_p).\;$ A K\"ahler metric $g$ (resp. a K\"ahler form $\omega$)
on a complex orbifold $M$  is a smooth metric on $M_{reg}$ such
that for any $p \in M,\; {\pi_p}^* g$ (resp.  K\"ahler form
${\pi_p}^* \omega$ )\footnote{Note ${\pi_p}^*$ is only defined
away from the fixed point set of $V_p.\;$ Since the fixed point
set is at least codimension 2 or higher, any metric defined on
non-fixed point set of $V_p$ has a unique smooth extension on
$V_p$ if such an extension exists. This definition essentially
says that a metric is smooth in the orbifold sense if such an
extension always exists in each uniformization system of the
underlying K\"ahler orbifold structure. } can extends to  a
smooth K\"ahler metric (resp.  smooth K\"ahler form) on $V_p.\;$
\end{defi}

\begin{defi} A function $f$ is called a smooth function on an
orbifold $M$  if for any $p \in M, f\circ \pi_p$ is a smooth
function on $V_p.\;$
\end{defi}
Similarly, one can define any tensor to be smooth on $M $  if its
pre-image on each local uniformization system  is smooth.
  Clearly, the curvature tensor and the Ricci
tensor of any smooth metric on orbifolds, as well as their
derivatives, are smooth tensors. A complex orbifold admits a
K\"ahler metric is called a K\"ahler orbifold.

\begin{defi} A curve $c(t)$ on K\"ahler orbifold $M $  is called geodesic if near any point $p$ on it,
  $c(t)\bigcap U_p $ can be lifted to  a geodesic on $V_p$
  and at least one preimage of $c(t)$  is smooth in $V_p.\;$
Here  $U_p$ is any open connected neighborhood of $p$ over which
$(V_p, G_p, \pi_p)$ is a unifomization system.
\end{defi}

Under this definition, we have
\begin{prop} Any  minimizing
geodesic between two  regular points never pass any singular
point of the K\"ahler orbifold.
\end{prop}
\begin{proof}  Otherwise, we can argue that the geodesic is not minimizing. Suppose that $p$ is a singular point and $p\in U_p$ is a small open
set which is uniformized by $(V_p, G_p, \pi_p)$ with an
equivariant metric $g$ on $V_p.\;$ Suppose that a portion of
geodesic lies inside of $U_p$ is $ c(t): [-\epsilon, \epsilon]$
such that $A= c(- \epsilon), B = c(\epsilon) \in U_p$ and $p =
c(0).\;$ Assume that this geodesic is parameterized by arc length.
Thus the distance between $A$ and $B$ is $2 \epsilon,$ while the
distance between $A$ (or $B$)  and $O$ is $\epsilon.\;$ Without
loss of generality, we may assume that $A$ and $B$ are regular
points. Suppose that $\pi^{-1}(p) = O;\; \pi^{-1}(A) = \{A_1,A_2,
\cdots, A_{l_1} \} $ and
 $\pi^{-1}(B) = \{ B_1, B_2,
\cdots, B_{l_2}\}.\;$ Note that $\{A_1,A_2,\cdots, A_{l_1}\}$ and
$\{B_1,B_2\cdots, B_{l_2}\}$ are on the sphere of radius
$\epsilon$ which centers at $O.\;$  If $G_p$ is non-trivial (then
the preimages of $A$ and $B$ are not unique, i.e., $l_1 > 1$ and
$l_2 > 1.$ ), then there is at least one pair of $A_i, B_j (1\leq
i\leq l_1,\; 1 \leq j \leq l_2)$ such that the distance  between
the two points is shorter than $2 \epsilon $ on $V_p$\footnote{In
any ball of radius $1$ on any metric space,  the maximum distance
between any two points in the ball is $2$ which is the diameter of
the unit ball.  Fixed a point in the ball, then the minimal
distance from that point to any set of points in the ball is
strictly less than $2, $ if that set contains two or more points.
}. Suppose this geodesic is $\tilde{C}.\;$ Then $\pi_p(\tilde{C})$
is a geodesic (which connects $A$ and $B$) whose length is shorter
than $2 \epsilon.\;$ Thus $c(t)$ is not a minimizing geodesic
between $A$ and $B$ since $\pi_p(A_i) = A$ and $\pi_p(B_j) = B.\;$
\end{proof}
\subsection{K\"ahler Ricci flow on K\"ahler-Einstein orbifolds}
A K\"ahler-Einstein orbifold metric  is a metric on orbifold such
that the Ricci curvature is proportional to the metric.  A
K\"ahler orbifold with a K\"ahler-Einstein metric is called a
K\"ahler-Einstein orbifold.

\begin{theo} Let $M$ be  any K\"ahler Einstein orbifold. If there is
another  K\"ahler metric in the same cohomology class which has
non-negative bisectional curvature and positive at least at one
point, then the K\"ahler-Ricci flow converges to a
K\"ahler-Einstein metric with positive bisetcional curvature.
\end{theo}
We want to generalize our proof of Theorem 1.1 to the orbifold
case. Note that the analysis for K\"ahler orbifolds is exactly the
same as that for  K\"ahler manifolds (cf. \cite{DingTian92}). We
want to show that this theorem can be proved exactly like Theorem
1.1. First,  we need to set up some notations. Following Section
2.1, we use the K\"ahler form $\omega$ as a smooth K\"ahler form
on the orbifold $M.\; $ Locally on $M_{reg}, $ it can be written
as
\[
\omega = \sqrt{-1} \displaystyle \sum_{i,j=1}^n\;g_{i
\overline{j}} d\,z^i\wedge d\,z^{\overline{j}},
\]
where $\{g_{i\overline {j}}\}$ is a positive definite Hermitian
matrix function.  Denote by $\cal B$  the set of all real valued
smooth functions on $M\;$ in the orbifold sense (cf. Definition
8.6).  Then the K\"ahler class $[\omega]$ consists of all K\"ahler
form which can be expressed as \[ \omega_{\varphi} = \omega + i
\partial \overline\partial \varphi
 > 0
\]
on $M$ for some $ \varphi \in \cal B.\;$ In other words, the space
of all K\"ahler potentials in this K\"ahler class is
\[
  {\cal H} = \{ \varphi \in  {\cal B}  \mid \omega_{\varphi} = \omega + i  \partial \overline\partial \varphi >
  0\}.
\]
  The Ricci form for $\omega $ is:
  \[
    Ric(\omega) = - \sqrt{-1} \partial \overline{\partial} \log \omega^n.
  \]
  As in the case of smooth manifolds,  $[\omega]$ is the canonical K\"ahler class if $\omega$ and the Ricci
  form $Ric(\omega)$ is in the same cohomology class after proper rescalling.  In the
  canonical K\"ahler class,  consider the K\"ahler Ricci flow

\[
   {{\partial \varphi} \over {\partial t }} =  \log {{\omega_{\varphi}}^n \over {\omega}^n } + \varphi - h_{\omega} ,
\]
where $h_{\omega}$ is defined as in Section 2.2. Clearly, this
flow preserves the structure of K\"ahler orbifold, in particular,
preserves the K\"ahler class $[\omega].\;$ Examining our proof of
Theorem 7.2, the following three parts are crucial
\begin{enumerate}
\item The preservation of positive bisectional curvature under the
K\"ahler Ricci flow.
\item The introduction of a set of new functionals $E_k$ and new
invariants $\Im_k (k=0,1,\cdots,n).\;$
\item The uniform estimate on the diameter; consequently, the
uniform control on the Sobolev constant and the Poincare constant.
\end{enumerate}

To extend these to the case of K\"ahler orbifolds, we really need
to make sure that the following tools for geometric analysis hold
in the orbifold case:

\begin{enumerate}
\item Maximum principle for smooth functions and tensors on
K\"ahler orbifold(cf. Definition 8.7). \item Integration by parts
for smooth functions/tensors in the orbifold case.
 \item The
second variation formula for any smooth geodesics(cf. Proposition
8.8).
\end{enumerate}

By our definition of K\"ahler orbifolds, it is not difficult to
see that the maximum principle holds  on orbifolds. Thus, Theorem
2.1 still holds in the orbifold case. In other words, the
bisectional curvature of the evolved metric is strictly positive
after the initial time,  if the initial metric has non-negative
bisectional curvature and positive at least at one point.
Moreover, the integration by parts on orbifold holds for any
smooth function on $M$  with smooth metrics in the orbifold sense.
Thus, our definitions of new functionals $E_0, E_1, \cdots E_n$
can be carried over to this K\"ahler orbifold setting without any
change.  Moreover, the formula for their derivatives still holds.
In particular,  $E_0$ and $E_1$ are decreasing strictly under the
K\"ahler Ricci flow. Furthermore, the set of invariants
$\Im_0,\Im_1,\cdots, \Im_n$ are well defined and vanish on any
K\"ahler-Einstein orbifold. Since Tian's inequality holds on any
K\"ahler Einstein orbifold, then Prop. 2.14, Corollary 2.16 hold
as well. Finally, the second variation formula for minimizing
geodesic between any two regular points on K\"ahler orbifolds is
exactly the same as the formula on smooth manifold(cf. Prop. 8.8).
Thus we can use the same set of ideas in Section 3 to estimate
diameter \footnote{In the proof of Lemma 3.5, without loss of
generality, may assume $p, q \in M_{reg}$ where the diameter $D=
d(p,q).\;$ Furthermore, we may assume $A_1 = B_{p,r} \subset
M_{reg}$ and $A_2 = B_{q,r} \subset M_{reg}.\;$ According to Lemma
8.8,  any minimizing geodesic between $A_1$ and $A_2$ belong to
$M_{reg}.\;$ Consequently, we can use Lemma 3.1 of J. Cheeger and
T. Coldings to conclude the diameter bounds as in the smooth case.
}; consequently, the Sobolev constant and the Poincare constant
can be uniformly controlled as well. The rest of arguments in our
proof of Theorem 7.2 can be extended to the orbifold case
directly. Thus we can prove Theorem 8.9 for K\"ahler-Einstein
orbifolds.

\subsection{K\"ahler Einstein orbifolds with constant positive bisectional
curvature} In this subsection, we want to prove the following

\begin{theo} Let $M$ be  any K\"ahler orbifold. If there is a
K\"ahler-Einstein metric with constant positive bisectional
curvature,  then it is a global quotient of $\CC P^n.\;$
\end{theo}
Suppose $\overline{g}$ is the standard Fubini-Study metric on $\CC
P^n$ with constant bisectional curvature. Suppose $g$ is a
K\"ahler-Einstein metric on $M$ with constant bisectional
curvature.  Normalize the bisectional curvature of $g$ on $M$ and
of $\overline{g}$ on $\CC P^n$  so that both bisectional
curvature is $1.\;$  Consequently, the conjugate radius of $\CC
P^n$ is $\pi.\;$ Let $p$ be any regular point in $M.\;$ By
definition, let $U_p$ be a small neighborhood of $p$ and $ (V_p,
G_p, \pi_p)$ be the uniformization system. Since $p \in M_{reg},
$  then $G_p$ is trivial group.  Consider $ g' = {\pi_p}^* g$ as a
K\"ahler metric with constant bisectional curvature on $V_p.\;$
If we choose $U_p$ sufficiently small, then $(V_p,g')$ is an open
subset of $(\CC P^n, \overline{g})$ with the induced  metric from
$(\CC P^n, \overline{g}).\;$  In the following, we will drop
notation $g'$ and use $\overline{g}$ only.  Our goal is to extend
$\pi_p$
into a local isometric map from $\CC P^n$ to $M.\;$\\

Next we set up some notations first.  Denote by $q$ the pre-image
of $p.\;$ Consider
\[\begin{array}{llcl}
   & \CC P^n &  & M \\
       &  \bigcup &  & \bigcup \\
   \pi_p:    &  V_p      & \rightarrowtail  & U_p.
   \end{array}
\]

Now, we want to lift this map $\pi_p$ to  a map from $\CC P^n $
to $ M.\;$  First, we need to rewrite this map in a  different
way:
\[
\begin{array} {llll} exp_q: & T_q(\CC P^n) & \rightarrow & \CC P^n \\
&             &             &   \bigcup \\
                     &             &             &   V_p \\
                     &    \downarrow  id   &             &   \downarrow  \pi_p  \\
                     &              &          &   U_p   \\
                     &             &             &   \bigcap \\
                  exp_p:   &  T_p(M)     & \rightarrow  &  M.

\end{array}
\]
Set \[ \Pi = exp_p \circ id \circ \exp_q^{-1}.\]
 Then at least
$\Pi$ is defined in $V_p,$ and
\begin{equation}
\Pi = \pi_p = exp_p \circ id \circ \exp_q^{-1},\qquad {\rm in}\;
V_p. \label{eq:extension}
\end{equation}
Consider the open ball of radius $\pi$ in $T_q(\CC P^n)$ which we
will denote it by $B_{\pi}.\;$ Then $\exp_q^{-1} $ is well defined
on $ exp_q(B_{\pi}) \subset \CC P^n.\;$ The image of $\partial
B_{\pi}$ under the exponential map is a projective subspace of
codimension $1, $ which will be denoted as $\CC
P^{n-1}_{\infty}.\;$ Then
\[
exp_q(B_{\pi}) = \CC P^n \setminus \CC P^{n-1}(\infty). \] We
claim that we can extend the map $\Pi$ in this way to $\CC P^n
\setminus \CC P^{n-1}(\infty)$ via  Formula (\ref{eq:extension}).
The key step is the following lemma (in the following arguments,
we abuse notation by using  letters $p$ and $q$ for generic
points on $M.\;$ ).

\begin{lem} Any smooth geodesic on $M$ can be extended uniquely and
indefinitely. In particular, it can be extended uniquely (before
the length $\pi$\footnote{We are interested in the unique
extension up to length $\pi$ since it is the conjugate radius of
any K\"ahler metric with constant bisectional curvature 1.}).
\end{lem}
\begin{proof} Suppose $c(t):[0,a]$ is a geodesic defined on $M$
with length $a > 0.\;$  If $ c(a) \in M_{reg},\;$  then it can
easily be extended as usual. If $c(a) \in U_p $ for some $p \in
M_{sing},$ in particular, if $c(a) \in M_{sing}$,  we want to
extend the geodesic uniquely as well. Consider the part of
geodesic $c(t) \bigcap U_p;$ And still denotes it as $c(t).\;$
Suppose that $U_p$ is uniformized by $(V_p, G_p, \pi_p).\;$ For
convenience, the pull backed metric $g'_p = {\pi}^* g$ is a smooth
metric on $V_p$ and $G_p$ acts isometrically on $(V_p, g'_p).\;$
Consider its pre-images $\tilde{c}(t)$ in $V_p$ under $\pi_p$
(note that $\pi_p$ is a local isometric map from $(V_p,g'_p)$ to
$(U_p, g),$ in particular, if we restrict the map to $M_{reg}
\bigcap U_p.$). Although the preimages are not unique in $V_p,\; $
each preimage $\tilde{c}(t)$ has a unique extension on $V_p.\;$
More importantly, the images of these geodesic extensions on $V_p$
under $\pi_p$  are unique in $U_p.\;$ Therefore, the geodesic
$c(t)$ is also extendable uniquely in this setting.
\end{proof}
In fact, we have the following
\begin{cor} Any geodesic in a K\"ahler orbifold with constant
bisectional curvature can be extended long enough to become a
closed geodesic. If the bisectional curvature is $1,$ then the
length of each closed geodesic is either $2 \pi$ or ${{2\pi}\over
l}$ for some integers $l.\;$  Moreover, there exists a maximum
integer of all such integers $l,\;$ denoted by $l_{max}.\;$ Then
the conjugate radius of the K\"ahler orbifold with constant
bisectional curvature $1$  is $ {\pi \over l_{max}}.\;$
\end{cor}
\begin{lem} $\Pi$ can be extended to be a global map from $\CC
P^n$ to $M.\; $ Moreover,  $\Pi$ is a local isometry  from an
open dense set $ (\Pi^{-1}(M_{reg}), \overline{g})$ in $ \CC P^n$
to $(M_{reg}, g).\;$
\end{lem}
 \begin{proof}
Consider the open ball $B_{\pi} \subset T_q(\CC P^n)$ of radius
$\pi.\;$ Then the closure of $exp_q(B_{\pi})$ is just the whole
$\CC P^n.\;$ By the preceding  lemma,  $\Pi$ can be defined in the
open set $exp_p (B_{\pi}).\;$ Next taking the closure of this
map, we define a map from $\CC P^n$ to $M.\;$ Moreover, $\Pi$ is
a local isometry from $(\Pi^{-1} (M_{reg}), \overline{g})$ to
$(M,g).\;$
\end{proof}

Now we want to prove the following lemma
\begin{lem} For any point $A$ on $M,$ there exists only a
finite number of preimages on $\CC P^n.\;$
\end{lem}
\begin{proof} Otherwise, there exists an infinite number of
pre-images of the point $A$ on $\CC P^n.\;$ This set of infinite
number points must have a concentration point on $\CC P^n.\;$ In
particular, for any small $\epsilon > 0,$ there are at least two
preimages of $A$ such that the distance between these two points
in $\CC P^n$ is less than $\epsilon.\;$ Consider the image of the
minimal geodesic which connects these two pre-image points on $\CC
P^n$ under $\Pi.\;$ We obtain  a geodesic loop centered at point
$A$ whose length is less than $\epsilon.\;$ This violates the fact
that the conjugate radius of $(M,g)$ is at least $ {\pi \over
l_{max}} $ (cf. Corollary 8.12).   Thus the lemma holds.
\end{proof}

\begin{lem} For any $p \in M, $ let $U_p$ be a small neighborhood
of $p$ and $(V_p, G_p, \pi_p)$ be the uniformization system over
$U_p.\;$  Let $W_p$ be any connected component of
$\Pi^{-1}(U_p).\;$Then there exists a finite group $G'_p$ acting
isometrically on $(W_p, \overline{g})$ such that $(W_p, G'_p,
\Pi\mid_{W_p}) $ is a uniformization system over $U_p,$ which is
equivalent to $(V_p, G_p, \pi_p)$ (In particular, if we choose a
different connected components of $\Pi^{-1}(U_p),$ then the
induced uniformization system $(W_p, G'_p, \Pi\mid_{W_p})$ is
invariant up to isometries.).
\end{lem}
\begin{proof} Set $f = \Pi\mid_W$ and $g'= {\pi}^* g.\;$
Then $g'$ is a smooth K\"ahler metric  with constant bisectional
curvature $1$ on $V_p,\;$ where $G_p$  acts isometrically on $V_p
$ with respect to this metric $g'.\;$ There exists a smooth
lifting of $f$ to $\tilde{f}: W_p \rightarrow V_p $ such that
$\pi_p \circ \tilde{f} = f.\;$ It is easy to verify that
$\tilde{f}$ is an isometric map from $(W_p, \overline{g})$ to
$(V_p,g').\;$ Moreover, it can be proved that $\tilde{f}$ is an
one-to-one map from $W_p$ to $V_p.\;$ Now consider the pull back
group $G' = \tilde{f}^{-1} G_p;$ and define the group action via
$\tilde{f}.\;$ Since $G_p$ acts isometrically on $(V_p, g'),$
then $ G'_p$ acts isometrically on $(W_p, \overline{g}).\;$ By
definition, $(W_p, G'_p, f)$ is a uniformization system over
$U_p$ which is equivalent to the original uniformization system
$(V_p, G_p, \pi_p).\;$
\end{proof}

\begin{lem}
For any $p\in M,$ let $U_p$ be a small neighborhood with
$(W_p,G_p,\Pi\mid_{W_p})$  a uniformization system where $(W_p,
\overline{g}) \subset (\CC P^n, \overline{g}).\;$ Then the fixed
point set of $G_p$ is a totally geodesic K\"ahler submanifold of
$(W_p, \overline{g}) \subset (\CC P^n, \overline{g}).\;$\end{lem}
\begin{proof} Following properties of isometric group actions on
K\"ahler manifold. 
\end{proof}
{\rm This lemma has an immediate}
\begin{cor} Consider $W_{sing} = \Pi^{-1} (M_{sing}).\;$ Then
$W_{sing}$ is a union of $CP^k$ for some $k \leq n-1.\;$  When
$k=0$, it is the preimage of isolated singular points on $M.\;$
\end{cor}
\begin{proof} The only totally geodesic K\"ahler submanifold in
$(\CC P^n, \overline{g})$ is $\CC P^k$ for some $k \leq n-1.\;$
\end{proof}
{\rm Denote by $Aut(\CC P^n)$  the holomorphic transformation
group of  $\CC P^n.\;$ Then we have}
\begin{lem}
For any $p\in M,$ let $U_p$ be a small neighborhood and  $(W_p,
G_p, \Pi\mid_{W_p}) $ be a uniformization system over $U_p,$ where
$(W_p, \overline{g}) \subset (\CC P^n, \overline{g}).\;$ For any
$\sigma \in G_p, $ it can be extended to be  a  group element in
$Aut (\CC P^n),$ and we still denote it as $\sigma.\;$ Moreover,
$\Pi \circ \sigma = \Pi$ on $\CC P^n.\;$
\end{lem}
\begin{proof} It is easy to see that $\sigma $ can be extended
uniquely to an element in $ Aut(\CC P^n)$ which acts isometrically
on $(\CC P^n, \overline{g}).\;$ Consider two local isometries
from $(\CC P^n, \overline{g})$ to $(M,g):\; \Pi$ and $\Pi \circ
\sigma.\;$ Since the two maps agree on an open set $W_p \subset
\CC P^n, $ they must agree on all $\CC P^n.\;$
\end{proof}

{\rm From now on, we may view $G_p$ as a subgroup of $Aut(\CC
P^n)$ directly.  Now we are ready to give a proof of Theorem
8.10.}
\begin{proof}
For any $p\in M,$ let $U_p$ be a small neighborhood with a
uniformization system $(W_p, G_p, \Pi\mid_{W_p}), $ then the
preceding lemma implies that $G_p$ is a subgroup of $Aut(\CC
P^n).\;$ If $p \in M_{reg}, $ then $G_p$ is trivial. If $p_1, p_2
\in M_{sing} $ and is near to each other, then $G_{p_1} =
G_{p_2}$ by continuity. Consequently, for any $p_1, p_2 \in
M_{sing}$ such that the fixed point sets of $W_{p_1}$ and
$W_{p_2}$ belong to the same connected component of
$\Pi^{-1}(M_{sing}),\;$ then $G_{p_1} = G_{p_2} \subset Aut(\CC
P^n). \;$ Consider $G \subset Aut(\CC P^n) $ to be the subgroup
generated by all such $G_p$'s. Then $G$ acts isometrically on
$(\CC P^n, \overline{g})$ and
\begin{equation}
   \Pi \circ \sigma = \Pi, \qquad \forall \sigma \in G \subset
   Aut(\CC P^n).
   \label{eq:globalcover}
 \end{equation}
This induces a covering map
\[
   \CC P^n / G \rightarrow M.
\]
By this explicit construction, one can verify directly that this
is an orbifold isomorphism. Consequently,  $M$ is a global
quotient of $\CC P^n$ by this group $G.\;$ The only thing left is
to show that $G$ is of finite order, which follows directly from
 equation (\ref{eq:globalcover}) and Lemma 8.14.
\end{proof}

\subsection{Pinching theorem for bisectional curvature}

{\rm In this subsection, we want to prove the following lemma}

\begin{lem}   If $g$ is a K\"ahler-Einstein metric with strictly
positive bisectional curvature on a K\"ahler orbifold, then $g$
has constant bisectional curvature.
\end{lem}
{\rm This lemma was first proved by Berger \cite{Berger65} on
K\"ahler manifolds. We note that his proof can be easily modified
for  K\"ahler orbifolds.  For reader's convenience, we include a
proof here.}

\begin{proof}  We begin with a
simple observation: In any K\"ahler orbifold, any
K\"ahler-Einstein metric satisfies the following elliptic
equation :
\[
\bigtriangleup R_{i \overline{j} k \overline{l}} + R_{i
\overline{j} p \overline{q}} R_{q \overline{p} k \overline{l}} -
R_{i \overline{p} k \overline{q}} R_{p \overline{j} q
\overline{l}} + R_{i \overline{l} p \overline{q}} R_{q
\overline{p} k \overline{j}} - R_{i \overline{j} k \overline{l}}
= 0.
\]

Define a new symmetric tensor $T_{i\overline{j} k \overline{l}}$
as
\[
  T_{i\overline{j} k \overline{l}} = g_{i \overline{j}} g_{k \overline{l}} +  g_{i \overline{l}}
\;g_{k \overline{j}}.
\]
And for any fixed $\epsilon \in (0,{1\over {n+1}})$, we define
\[
  Q_{i \bar j k \bar l} = R_{i \bar j k \bar l} - \epsilon T_{ i
  \bar j k \bar l}.
\]

Note that $T_{i \bar j k \bar l}$ is parallel in the manifold. By
a direct but tedious calculation,  we arrive at the following

\begin{eqnarray}
 & & - \triangle Q_{i \bar j k \bar l}  \nonumber \\
 & = &  Q_{i \overline{j} p \overline{q}} Q_{q \overline{p} k \bar
l} +Q_{i \overline{l} p \overline{q}} Q_{q \overline{p} k \bar j}
- Q_{i \overline{p} k \overline{q}} Q_{p \overline{j} q
\overline{l}} \nonumber \\ & & \qquad \qquad + \epsilon ( 1 -
(n+1) \epsilon) T_{i \bar j k \bar l} - Q_{i \overline{j} k
\overline{l}}. \label{iden5}
\end{eqnarray}
Suppose that the bisectional curvature of $g$ is not constant.
Note that for $\epsilon = {1\over {n+1}},$ we have
\[
g^{i \bar j} Q_{i \bar j k \bar l} = g^{k \bar l} Q_{i \bar j k
\bar l} = 0.
\]
Thus if $R_{i \bar j k \bar l} > 0,$ there exists a small
positive $\epsilon \in (0, {1\over {n+1}} ),$ such that $Q_{i \bar
j k \bar l} \geq 0$ in the whole manifold and vanishes in some
direction at some points. In other words, there exists a point
$x_0 \in M$ and two vectors $v_0, w_0 \in T_{x_0} M$ such that
\[  Q_{i \bar j k \bar l}(x_0) v_0^{\bar i} v_0^j w_0^{\bar k}
w_0^l = 0
\]
and for any other point $x\in M$ and any other pair of vectors
$v, w \in T_x M,$ we have
\[
  Q_{i \bar j k \bar l}(x) v^{\bar i} v^j w^{\bar k}
w^l \geq 0.
\]
Now consider a pair of parallel vector fields $v,w$ in a small
neighborhood of $x_0$ such that
\[
   v^i_{,j} = v^i_{,\bar j} = w^i_{,j} = w^i_{,\bar j}=0,
\]
where
\[
  v = \displaystyle \sum_{i=1}^n\; v^i\; {{\p }\over {\p
  z^i}},\qquad {\rm and}\qquad  w = \displaystyle \sum_{i=1}^n\; w^i\; {{\p }\over {\p
  z^i}}.
\]
Furthermore, $v(x_0) = v_0$ and $w(x_0) = w_0.\;$  Consider the
scalar function
\[
  Q =  Q_{i \bar j k \bar l}(x) v^{\bar i} v^j w^{\bar k}
w^l
\]
in a neighborhood of $x_0.\;$ Clearly, $Q$ achieves minimum in
$x_0.\;$ Thus the maximum principle implies that
\[
\begin{array}{lcl}
- \triangle Q & = &  \left(\triangle Q_{i \bar j k \bar l} v^{\bar
i} v^j w^{\bar k} w^l\right)_{x=x_0} \\
& = & \left( \triangle Q_{i \bar j k \bar l}  \right)\mid_{x=x_0}
v^{\bar i} v^j w^{\bar k} w^l \leq 0 \end{array}.
\]
Plugging this into the equation (\ref{iden5}), we obtain
 \[
\begin{array}{lcl}
& & - \triangle Q_{i \bar j k \bar l} v^{\bar i} v^j w^{\bar k}
w^l\\ & = & Q_{i \overline{j} p \overline{q}} Q_{q \overline{p} k
\bar l}  v^{\bar i} v^j w^{\bar k} w^l  + Q_{i \overline{l} p
\overline{q}} Q_{q \overline{p} k \bar j}  v^{\bar i} v^j w^{\bar
k} w^l  - Q_{i \overline{p} k \overline{q}} Q_{p \overline{j} q
\overline{l}}  v^{\bar i} v^j w^{\bar k} w^l \\
& & \qquad \qquad  + \epsilon ( 1 - (n+1) \epsilon) T_{i \bar j k
\bar l} v^{\bar i} v^j w^{\bar k} w^l - Q_{i \overline{j} k
\overline{l}}  v^{\bar i} v^j w^{\bar k} w^l.
\end{array}
\]
Define the following linear operator at point $x_0\;$
\[
  A_{i \bar j} = R_{i \bar j k \bar l} v^{\bar k} v^l,\qquad {\rm
  and}\qquad C_{i \bar j} =  R_{k \bar l i \bar j } w^{\bar k}
  w^l
 \]
 and
 \[
 M_{i \bar j} = R_{ i \bar p q \bar j} v^p w^{\bar q},\qquad
 {\rm and}\qquad  M_{i  j} = R_{ i \bar p j \bar q} v^p w^{ q}.
 \]
 Plugging these  into the above equation and evaluate at  $x_0, $
 we have
 \[\begin{array}{lcl}
0 & \geq & - \triangle Q_{i \bar j k \bar l}(x) v^{\bar i} v^j
w^{\bar k} w^l \mid_{x=x_0} \\ & = &  A_{p \bar q} C_{q \bar p} +
M_{p \bar q} M_{q \bar p} - M_{p q}M_{\bar p \bar q} + \epsilon
(1 - (n+1)\epsilon) |v|^2 |w|^2.
\end{array}
 \]
 By a calculation of N. Mok, we have
 \[
 A_{p \bar q} C_{q \bar p} \geq M_{p \bar q}
M_{q \bar p} + M_{p q}M_{\bar p \bar q}.
 \]
 Since $0 < \epsilon < {1\over {n+1}},\;$ then
 \[
0  \geq  - \triangle Q_{i \bar j k \bar l}(x_0) v^{\bar i} v^j
w^{\bar k} w^l \geq  \epsilon (1 - (n+1)\epsilon) |v|^2 |w|^2 > 0.
 \]
 This is a contradiction!  Thus, $\epsilon = {1\over {n+1}}.\;$ In
 this case,
\[
g^{i \bar j} Q_{i \bar j k \bar l} = g^{k \bar l} Q_{i \bar j k
\bar l} = 0
\]
and \[ Q_{i \bar j k \bar l}  \geq 0. \]
 Thus
\[
Q_{i \bar j k \bar l} (x_0) = 0.
\]

Since $x_0$ is the minimum for $R_{i \bar j k\bar l},$ we obtain
that
\[
  R_{i \bar j k \bar l} \equiv {1\over {n+1}} T_{i \bar j k \bar l}\equiv {1\over {n+1}} \left(g_{i \bar j} g_{ k \bar
  l} + g_{i \bar l} g_{ k \bar j}\right).
\]
\end{proof}

{\rm From the proof, we can actually prove slightly more:}
\begin{lem} If the K\"ahler Ricci flow converges to a unique K\"ahler
Einstein metric exponentially fast (Prop. 8.3), and if the
bisectional curvature remains positive before taking limit, then
the limit K\"ahler Einstein metric has positive bisectional
curvature.
\end{lem}

{\rm Combining this lemma, Lemma 8.19 and Theorem 7.3, we can
prove}
\begin{theo} For any integer $l>0,\;{{\partial \varphi} \over {\partial t}} $
   converges exponentially fast to $0\;$ in any $C^l$ norm. Furthermore, the
   K\"ahler Ricci flow converges exponentially fast to a unique
   K\"ahler Einstein metric with constant bisectional curvature on any K\"ahler-Einstein manifolds.
\end{theo}

\section{Concluding Remarks}
{\rm In this  section, we want to prove our main Theorem 1.1,
Corollary 1.3 and Theorem 1.4. Note that the proof of Theorem 1.1
and Corollary 1.3 is similar to the proof we gave in the final
section of our earlier paper \cite{chentian001} except in the
final step, where we need to use Berger's theorem (Lemma 8.19) and
Theorem 8.21 to show that any K\"ahler Einstein metric with
positive bisectional curvature must be a space form. We will skip
this part and just give a proof for Theorem 1.4.}

\begin{proof}   For any K\"ahler metric  in the
canonical K\"ahler class  such that it has non-negative
bisectional curvature on $M$(K\"ahler-Einstein orbifold) but
positive bisectional curvature at least at one point, we apply the
K\"ahler Ricci flow to this metric on $M.\;$ Theorem 2.1 still
holds in the orbifold case.  In other words, the bisectional
curvature of the evolved metric is strictly positive over all the
time. By our Theorem 8.9 and Lemma 8.20, the K\"ahler Ricci flow
converges exponentially to a unique K\"ahler-Einstein metric of
positive bisectional curvature. According to Lemma 8.19, any
K\"ahler-Einstein metric of positive bisectional curvature on a
K\"ahler orbifold must have a constant positive bisectional
curvature.  Moreover, using Theorem 8.10, we arrive at the
conclusion that $M$ must be
a global quotient of $\CC P^n $ by a finite group action. \\

Furthermore, this also proves that any K\"ahler metric with
nonnegative bisectional curvature on $M$ and positive at least at
one point is path connected to a K\"ahler-Einstein metric of
constant positive bisectional curvature. Note that all the
K\"ahler-Einstein metrics are path connected by automorphisms
\cite{Ma85}. Therefore, the space of all K\"ahler metrics with
nonnegative bisectional curvature on $M$ and positive at least at
one point, is path connected. Similarly, using Theorem 2.2 and
our Theorem 1.4, we can show that all of K\"ahler metrics with
nonnegative curvature operator on $M$ and positive at least at
one point is path connected. Note that the nonnegative curvature
operator implies the nonnegative bisectional curvature.
\end{proof}

\begin{rem} Combining our main Theorem 1.1 and Theorems 2.1, 2.2, we
can easily generalize Corollary 1.2 to the case that the
bisectional curvature (or curvature operator) is only assumed to
be non-negative. We can show the flow converges exponentially
fast to a unique K\"ahler-Einstein metric with non-negative
bisectional curvature. Then, following an earlier work of Zhong
and Mok, the underlying manifold must be compact symmetric
homogeneous manifold.
\end{rem}

{\rm Next we want to propose some future problems. Some of them
may not be hard to solve.}

\begin{q} It is clear that $E_1$ plays a critically important
role in proving the convergence theorem of this paper. Note that
the Ricci flow is a gradient-like of the functional $E_1.\;$ It
will be interesting to study the gradient flow of $E_1.\;$ \end{q}

{\rm  Consider the expansion formula in $t$:
$$ \left ( \omega_\varphi + t \;{\rm Ric} (\omega_\varphi) \right
)^n = \left( \sum_{k=0}^n \sigma_k(\omega_\varphi) t^k \right )
\omega_\varphi^n.$$ Clearly, $\sigma_0(\omega_\varphi)=1$,
$\sigma_1(\omega_\varphi) = R(\omega_\varphi)$, the scalar
curvature of $\omega_\varphi$. The equation for the gradient flow
of $E_1$ is

\begin{equation}
{{\p \varphi(t)}\over {\p t}} =2  \Delta_\varphi R(\omega_\varphi)
-(n-1)\sigma_{2}(\omega_\varphi ) - c_1.
\end{equation}
Here $c_1$ is some constant which depends only on the K\"ahler
class. Clearly, this is a 6 order parabolic equation.}

\begin{q} As Remark 1.5 indicates, what we really need is the
positivity of Ricci curvature along the K\"ahler Ricci flow.
However, it is not expected that the positivity of Ricci
curvature is preserved under the K\"ahler Ricci flow except on
Riemann surfaces. The positivity of bisectional curvature is a
technical assumption to assure the positivity of Ricci curvature.
It is very interesting to extend Theorem 1.1 to metrics without
the assumption on the bisectional curvature.
\end{q}

\begin{conj} (Hamilton-Tian) On a K\"ahler-Einstein manifold, the K\"ahler-Ricci
flow converges to a K\"ahler-Einstein metric. On a general
K\"ahler manifold with positive first Chern class, the K\"ahler
Ricci flow will converge, at least by taking sequences, to a
K\"ahler Ricci soliton modulo diffeomorphism. Note that the
complex structure may change in the limit and the limit may have
mild singularities.
\end{conj}

\begin{q} Is the positivity of the sectional curvature preserved
under the K\"ahler-Ricci flow?
\end{q}

\begin{q} For any holomorphic vector field and K\"ahler class
$[\omega],$ are the invariants ${\cal{I}}_k (X, [\omega])$
independent? In the non-canonical class, we expect these
invariants to be different. Note that one can derive localization
formulas for  these invariants as what are done  in \cite{tian98}.
\end{q}

\begin{q} According to our Theorem 1.4,  any K\"ahler-Einstein orbifold
with positive bisectional curvature is necessarily biholomorphic
to a global quotient of $\CC P^n.$ What happens if we drop the
K\"ahler-Einstein condition?
\end{q}


\begin{thebibliography}{10}

\bibitem{Bando84}
S.~Bando.
\newblock On the three dimensional compact {K}\"ahler manifolds of nonnegative
  bisectional curvature.
\newblock {\em J. {D}. {G}.}, 19:283--297, 1984.

\bibitem{Ma85}
S.~Bando and T.~Mabuchi.
\newblock Uniqueness of {E}instein {K}aehler metrics modulo connected group
  actions.
\newblock {\em Algebraic geometry, Sendai, 1985, Adv. Stud. Pure Math.},
  10:11--40, 1987.

\bibitem{Berger65}
M.~Berger.
\newblock {\it Sur les vari{\'e}t{\'e}s d'{E}instein compactes}.
\newblock {\em C. R. IIIe R{\'e}union Math. Expression latine, Namur},
  35--55, 1965.

\bibitem{Cao85}
H.~D. Cao.
\newblock Deformation of {K}\"ahler metrics to {K}\"ahler-{E}instein metrics on
  compact {K}\"ahler manifolds.
\newblock {\em Invent. Math.}, 81:359--372, 1985.

\bibitem{caoprivate}
H.~D. Cao and R.~Hamilton.
\newblock On the preservation of positive orthogonal bisectional curvature
  under the k\"ahler ricci flow, 1998.
\newblock private communication.

\bibitem{CC96}
J.~Cheeger and T.~H. Colding.
\newblock Lower bounds on{R}icci curvature and almost rigidity of wraped
  products.
\newblock {\em Ann. Math.}, 144:189--237, 1996.

\bibitem{chen994}
X.~X. Chen.
\newblock Calabi flow in {R}iemann surface revisited: a new point of views.
\newblock (6):276--297, 2001.
\newblock ``International {M}athematics {R}esearch {N}otices".

\bibitem{chentian001}
X.~X. Chen and G.~Tian.
\newblock Ricci flow on complex surfaces, 2000.
\newblock preprint,.

\bibitem{Chow91}
B.~Chow.
\newblock The {R}icci flow on the 2-sphere.
\newblock {\em J. {D}iff. {G}eom.}, 33:325--334, 1991.

\bibitem{croke80}
C.~Croke.
\newblock Some {I}soperimetric {I}nequalities and {C}onsequences.
\newblock {\em Ann. Sci. E. N. S., Paris}, 13:419--435, 1980.

\bibitem{DingTian92}
W.~Ding and G.~Tian.
\newblock K\"ahler-{E}istein metric and the generalized {F}utaki invariant.
\newblock {\em Inventiones mathematicae}, 110:315--335, 1992.

\bibitem{futaki83}
A.~Futaki.
\newblock An obstruction to the existence of {E}instein {K}\"ahler metrics.
\newblock {\em Inv. Math. Fasc.}, 73(3):437--443, 1983.

\bibitem{Hamilton82}
R.~Hamilton.
\newblock Three-manifolds with positive {R}icci curvature.
\newblock {\em J. {D}iff. {G}eom.}, 17:255--306, 1982.

\bibitem{Hamilton86}
R.~Hamilton.
\newblock Four-manifolds with positive curvature operator.
\newblock {\em J. {D}iff. {G}eom.}, 24:153--179, 1986.

\bibitem{Hamilton88}
R.~Hamilton.
\newblock The {R}icci flow on surfaces.
\newblock {\em Contemporary {M}athematics}, 71:237--261, 1988.

\bibitem{Hamilton93}
R.~Hamilton.
\newblock {\em The formation of singularities in the Ricci flow}, volume~II.
\newblock Internat. Press, 1993.

\bibitem{pli80}
P.~Li.
\newblock On the sobolev constant and the $p-$ spectrum of a compact
  {R}iemannian manifold.
\newblock {\em Ann. Sci. E. N. S., Paris}, 13:451--468, 1980.

\bibitem{liyau80}
P.~Li and S.~T. Yau.
\newblock Estimates of eigenvalues of a compact riemannian manifold.
\newblock In {\em Proceedings of {S}ymposia in {P}ure {M}athematics},
  volume~36, pages 205--239, 1979.

\bibitem{Marg86}
C.~Margein.
\newblock Positive {P}inched manifolds are space forms.
\newblock {\em Proceedings of Symp. in Pure Math.}, 44:307--328, 1986.

\bibitem{Mok88}
N.~Mok.
\newblock The uniformization theorem for compact {K}\"ahler manifolds of
  non-negative holomorphic bisectional curvature.
\newblock {\em J. {D}ifferential {G}eom.}, 27:179--214, 1988.

\bibitem{Mori79}
S.~Mori.
\newblock Projective manifolds with ample tangent bundles.
\newblock {\em Ann. of. {M}ath.}, 76(2):213--234, 1979.

\bibitem{Nish86}
S.~Nishikawa.
\newblock Deformation of {R}iemannian metrics and manifolds with bounded
  curvature ratios.
\newblock {\em Proceedings of Symp. in Pure Math.}, 44:345--352, 1986.

\bibitem{Yuan00}
Y.~B. Ruan.
\newblock String geometry and topology of orbifolds, 2000.
\newblock preprint, AG/0011149.

\bibitem{Siuy80}
Y.T. Siu and S.~T. Yau.
\newblock Compact {K}\"ahler manifolds of positive bisectional curvature.
\newblock {\em Invent. Math.}, 59:189--204, 1980.

\bibitem{spro001}
C.~Sprouse.
\newblock Integral curvature bounds and bounded diameter.
\newblock {\em Communication of {A}nalysis and {G}eometry}, 8(3):531--543,
  2000.

\bibitem{Tian97}
G.~Tian.
\newblock K\"ahler-{E}instein metrics with positive scalar curvature.
\newblock {\em Invent. Math.}, 130:1--39, 1997.

\bibitem{tian98}
G.~Tian.
\newblock Some aspects of {K}\"ahler {G}eometry, 1997.
\newblock Lecture note taken by {M}eike {A}keveld.

\bibitem{Yau78}
S.~T. Yau.
\newblock On the {R}icci curvature of a compact {K}\"ahler manifold and the
  complex {M}onge-{A}mpere equation, ${I}^*$.
\newblock {\em Comm. Pure Appl. Math.,}, 31:339--441, 1978.

\end{thebibliography}

\end{document}